\documentclass{birkmult}

\usepackage{amssymb,amsmath,amscd,graphicx,amsthm}

\newtheorem{theorem}{Theorem}
\newtheorem{lemma}[theorem]{Lemma}
\newtheorem{proposition}[theorem]{Proposition}

\theoremstyle{definition}

\theoremstyle{remark}
\newtheorem{remark}[theorem]{Remark}
\newtheorem{example}[theorem]{Example}
\numberwithin{equation}{section}
\numberwithin{theorem}{section}

\def\A{{\mathcal A}}
\def\AA{{\mathfrak A}}

\def\CC{{\mathcal E}}
\def\EE{{\mathcal E}}

\def\FF{{\mathcal F}}
\def\FFF{{\mathfrak F}}
\def\G{{\mathcal G}}
\def\H{\mathcal H}

\def\L{{\mathcal L}}

\def\P{{\mathcal P}}
\def\PP{{\mathfrak P}}
\def\R{{\mathbb R}}
\def\RR{{\mathcal R}}
\def\RS{{(\R\setminus 0)}}

\def\Z{{\mathbb Z}}

\def\g{\mathfrak g}

\def\ls{\le}

\def\one{\mathbf 1}

\def\wB{{\widetilde{B}}}
\def\wM{{\widetilde{M}}}
\def\wN{{\widetilde{N}}}

\def\wX{{\bar X}}
\def\hM{{\widehat{M}}}
\def\hN{{\widehat{N}}}

\def\hP{{\widehat{P}}}
\def\wx{{\widetilde{\bf x}}}
\def\x{{\bf x}}

\def\End{\operatorname{End}}

\def\Mat{\operatorname{Mat}}
\def\Net{\operatorname{Net}}

\def\Poi{{\{\cdot,\cdot\}}}

\def\Tr{\operatorname{Tr}}

\def\diag{\operatorname{diag}}

\def\grad{\mbox{grad}}

\def\rank{\operatorname{rank}}
\def\s{\operatorname{sign}}

\def\:{{:\ }}

\begin{document}

\title{ Poisson Geometry of Directed Networks in a Disk}

\author{Michael Gekhtman}

\address{Department of Mathematics, University of Notre Dame, Notre Dame,
IN 46556}
\email{mgekhtma@nd.edu}

\author{Michael Shapiro}
\address{Department of Mathematics, Michigan State University, East Lansing,
MI 48823}
\email{mshapiro@math.msu.edu}
 
\author{Alek Vainshtein}
\address{Department of Mathematics \& Department of Computer Science, University of Haifa, Haifa,
Mount Carmel 31905, Israel}
\email{alek@cs.haifa.ac.il}

\begin{abstract} We investigate Poisson properties of Postnikov's map from the space of  edge weights 
of a planar directed network into the Grassmannian. We show that this map is Poisson if the space of  edge weights is equipped with a representative of a 6-parameter family of universal quadratic Poisson brackets and the Grasmannian is viewed
as a Poisson homogeneous space of the general linear group equipped with an appropriately chosen
R-matrix Poisson-Lie structure. We also prove that Poisson brackets on the Grassmannian arising in this way
are compatible with the natural cluster algebra structure.

\end{abstract}


\maketitle

\bigskip 

\section{Introduction}

Directed planar graphs with weighted edges have been widely used
in the study of totally nonnegative matrices (\cite{KarlinMacGregor, Brenti, BFZ}). ( Reviews of the area can be found in \cite{FZ_Intel, Fallat}. )
In particular, a special kind of such graphs is a convenient tool
for visualizing in the $GL_n$ case Lusztig type  parametrizations of double Bruhat cells 
\cite{Lusztig, FZBruhat}. Each parametrization of this kind is  obtained via a factorization of an element of a cell into a product of elementary factors. The standard Poisson-Lie structure on reductive Lie groups is induced by a very simple Poisson bracket on factorization 
parameters \cite{Reshetikhin&Co,R2}. 

Recently, Postnikov used 
weighted directed planar graphs to parametrize cells in Grassmannians \cite{Postnikov}.
The main goal of this paper is to investigate Poisson properties of Postnikov's parametrization. Using a graphical interpretation of the Poisson-Lie property inspired by the Lusztig parametrization, we introduce   a natural family of Poisson brackets on the space of edge weights. These brackets induce a two-parameter family of Poisson brackets on the Grassmannian. Every Poisson bracket in this family is compatible (in the sense of \cite{GSV1,GSV2}) with the cluster algebra on the Grassmannian described in \cite{GSV1,Scott} and, on the other hand, endows
the Grassmannian with a structure of a Poisson homogeneous space with respect
to the natural action of the general linear group equipped with an R-matrix Poisson-Lie
structure.

This paper is the first in the series devoted to the study of Poisson properties of weighted directed graphs on surfaces. The second paper \cite{GSV3} generalizes results of the current one
to graphs in an annulus. In this case, the analogue of Postnikov's construction leads
to a map into the space of loops in the Grassmannian. Natural Poisson brackets on
edge weights in this case  are  intimately connected to trigonometric R-matrix brackets
on matrix-valued rational functions. The third paper \cite{GSV4} utilizes particular graphs in an annulus to introduce a cluster algebra structure related to the coordinate ring of the space
of normalized rational functions in one variable. This space is birationally equivalent, via the Moser map \cite{moser},  to any minimal irreducible coadjoint orbit of the group of upper triangular  matrices associated with a Coxeter element of the permutation group. In this case, the Poisson bracket compatible with the cluster algebra structure coincides with the quadratic Poisson bracket studied in \cite{FayGekh1, FayGekh2} in the context of Toda flows on minimal orbits. We show that cluster transformations serve as B\"acklund-Darboux transformations between different minimal Toda flows. 
The fourth paper \cite{GSV5} solves, in the case of graphs in an annulus with one
source and one sink, the inverse problem of restoring the weights from the image
of the generalized Postnikov map. In the case of arbitrary planar graphs in a disk, this
problem was completely solved by Postnikov \cite{Postnikov} who proved that for a fixed minimal graph, the space of weights modulo gauge action is birational to its image and described all minimal graphs.
 To the contrary, already for simplest graphs in an annulus, the corresponding map  can only be shown to be finite. In \cite{Pavlo}, the inverse problem for the totally nonnegative matrices
is solved for a particular type of minimal graphs.
In \cite{GSV5} we describe all minimal graphs with one source and one sink.

 The paper is organized as follows.

In Section \ref{PPN}, we introduce a notion of a perfect  network in a disk and
associate with every such network a {\em matrix of boundary measurements}. 
Each boundary measurement is shown to be a rational function in edge weights   and admits a subtraction-free rational expression, see Proposition~\ref{sfree}.
In  Section \ref{PSmat}, we characterize all  {\it universal\/} Poisson brackets on the
space of edge weights   of a given network that respect the natural operation of concatenation of
networks, see Proposition \ref{6param}. Furthermore, we establish that the family
of universal brackets induces a linear two-parameter family of Poisson brackets on
boundary measurement matrices, see Theorem~\ref{PSM}. This family
depends on a mutual location of sources and sinks, but not on the network itself.
We  provide an explicit description of this family in Theorem~\ref{PSME}.
In Section \ref{PSongr}, we start by showing that if the sources and the sinks are not intermixed along the 
boundary of the disk, one can recover a 2-parametric family of R-matrices and the
corresponding R-matrix brackets on $GL_n$, see Theorem~\ref{PSGL}.
Next, the  boundary measurement map  defined by a network with $k$ sources and $n-k$ sinks is extended to the Grassmannian  boundary measurement map into the
Grassmannian $G_k(n)$.  
The Poisson family on boundary measurement matrices allows us to
equip the Grassmannian with a  two-parameter family of Poisson brackets
$\P_{\alpha,\beta}$ in such a way that for any choice of a universal Poisson 
bracket on edge weights there is a unique member of $\P_{\alpha,\beta}$ 
that makes the Grassmannian boundary measurement map Poisson, see Theorem~\ref{PSGr}. 
This latter family depends only on the number of sources and sinks 
and does not depend on their mutual location.
Finally, we give an interpretation of the natural  $GL_n$ action on  $G_k(n)$ in terms of
networks and establish that  every member of the Poisson family $\P_{\alpha,\beta}$
on  $G_k(n)$ makes $G_k(n)$ into a Poisson homogeneous space of $GL_n$ equipped
with the above described R-matrix bracket, see Theorem~\ref{poihomo}. In Section \ref{compclust},
we review the construction of the cluster algebra structure on an open cell in the Grassmannian given
in \cite{GSV1}. We then introduce face weights and use them to show that  every member of the Poisson family $\P_{\alpha,\beta}$ is compatible with this cluster algebra structure, see Theorem~\ref{thm:PoissComp}.

\vskip 1cm

\section{Perfect planar networks and boundary measurements}
\label{PPN}
\subsection{Networks, paths and weights}
\label{PandW}

Let $G=(V,E)$ be a directed planar graph with no loops and parallel edges drawn inside a disk (and
considered up to an isotopy) with the vertex set $V$ and the edge
set $E$. 
Exactly $n$ of its vertices are located on the boundary circle of
the disk. They are labelled counterclockwise $b_1,\dots,b_n$ and
called {\it boundary vertices\/}; occasionally we will write $b_0$ for $b_n$ and $b_{n+1}$ for $b_1$. Each boundary vertex is marked 
as a source or a sink. A {\it source\/} is
a vertex with exactly one outcoming edge and no incoming edges.
{\it Sinks\/} are defined in the same way, with the direction of the single edge
reversed. The number of sources is denoted by $k$, and the corresponding set of indices, by
$I\subset [1,n]$; the set of the remaining $m=n-k$ indices is denoted by $J$.
All the internal vertices of $G$ have degree~$3$ and are of two types: either they have exactly one
incoming edge, or exactly one outcoming edge. The vertices of the first type are called
(and shown on figures) {\it white}, those of the second type, {\it black}.

Let $x_1,\dots,x_d$ be independent variables.
A {\it perfect planar network\/} $N=(G,w)$ is obtained from $G$ as above by assigning a weight $w_e\in \Z(x_1,\dots,x_d)$
to each edge $e\in E$. In what follows we occasionally write ``network'' instead of ``perfect planar network''. 
Each network defines a rational map $w: \R^d\to \R^{|E|}$;
the {\it space of edge weights\/} $\EE_N$ is defined as the intersection of the image of $(\R\setminus 0)^d$ under $w$ with $(\R\setminus 0)^{|E|}$. In other words, a point in $\EE_N$ is a graph $G$ as above with edges weighted by nonzero reals obtained by specializing the variables $x_1,\dots,x_d$ in the expressions for $w_e$ to nonzero values.


\begin{figure}[ht]
\begin{center}
\includegraphics[height=4.5cm]{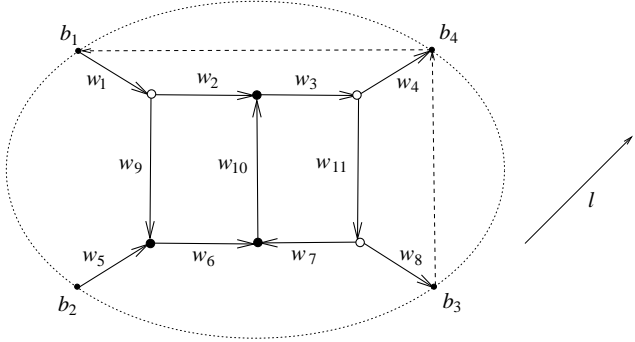}
\caption{A perfect planar network in a disk}
\label{fig:pergraph}
\end{center}
\end{figure}

An example of a perfect planar network is shown 
in Fig.~\ref{fig:pergraph}. It has two sources: $b_1$ and $b_2$,
and two sinks $b_3$ and $b_4$. Each edge $e_i$ is labelled by its
weight. The weights depend on four independent variables $x_1, x_2, x_3, x_4$ and are given by
\begin{align*}
&w_1=x_1^2/(x_2+1),\qquad w_2=x_2,\qquad w_3=x_2+1,\qquad w_4=x_1+x_3,\\
&w_5=x_3,\qquad w_6=x_3,\qquad w_7=x_3,\qquad w_8=x_4,\\
&w_9=1, \qquad w_{10}=1, \qquad w_{11}=1.
\end{align*}
The space of edge weights is the 4-dimensional subvariety in $(\R\setminus 0)^{11}$ given by equations $w_1w_3=(w_4-w_5)^2$, $w_3=w_2+1$, $w_5=w_6=w_7$, $w_9=w_{10}=w_{11}=1$ and condition $w_3\ne 1$.

A {\it path\/} $P$ in $N$ is an alternating sequence $(v_1, e_1,v_2,\dots,e_r, v_{r+1})$ of
vertices and edges such that $e_i=(v_i,v_{i+1})$ for any $i\in [1,r]$.
Sometimes we omit the names of the vertices and write $P=(e_1,\dots,e_r)$. 
A path is called a {\it cycle\/} if $v_{r+1}=v_1$ and a {\it simple cycle\/} if additionally $v_i\ne v_j$ for any other pair $i\ne j$.

To define the weights of the paths we need the following construction.
Consider a closed oriented polygonal plane curve $C$. Let $e'$ and $e''$ be two consequent oriented segments of $C$, and let $v$ be their common vertex. We assume for simplicity that for any such pair $(e',e'')$, the cone spanned by $e'$ and $e''$  is not a line; in other words, if $e'$ and $e''$ are collinear, then they have the same direction. Observe that since $C$ is not necessary simple, there might be other edges of $C$ incident to $v$ (see Figure~\ref{fig:concord} below). Let $l$ be an arbitrary oriented line. Define $c_l(e',e'')\in \Z/2\Z$ in the following way:
$c_l(e',e'')=1$ if the directing vector of $l$ belongs to the interior of the cone spanned by $e'$ and $e''$ and $c_l(e',e'')=0$  otherwise (see Figure~\ref{fig:concord} for examples). Define $c_l(C)$ as the sum of $c_l(e',e'')$ over all pairs of consequent segments in $C$. 
It follows immediately from Theorem~1 in \cite{GrSh} that $c_l(C)$ does not depend on $l$, provided $l$ is not collinear to any of the segments in $C$.
The common value of $c_l(C)$ for different choices of $l$ is denoted by $c(C)$ and called the {\it concordance number\/} of $C$. In fact, $c(C)$ equals $\mod 2$ the rotation number of $C$; the definition of the latter is similar, but more complicated.

\begin{figure}[ht]
\begin{center}
\includegraphics[height=4.5cm]{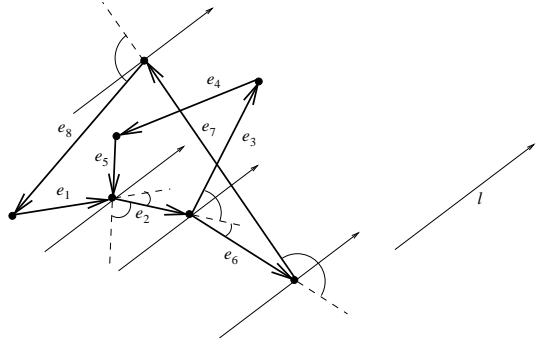}
\caption{To the definition of the concordance number: $c_l(e_1,e_2)=c_l(e_5,e_2)=0$; $c_l(e_2,e_3)=1, c_l(e_2,e_6)=0$; $c_l(e_6,e_7)=1, c_l(e_7,e_8)=0$}
\label{fig:concord}
\end{center}
\end{figure}

In what follows we assume without loss of generality that $N$ is drawn in such a way that all its edges are straight line segments 
and all internal vertices belong to the interior of the convex hull of the boundary vertices. 
Given a path $P$ between a source $b_i$ and a sink $b_j$, we define a closed polygonal curve $C_P$ by adding to $P$ the path between $b_j$ and $b_i$ that goes counterclockwise along the boundary of the convex hull of all the boundary vertices of $N$. Finally the weight of $P$ is defined as
\[
w_P=(-1)^{c(C_P)-1}\prod_{e\in P} w_e.
\]
The weight of an arbitrary cycle in $N$ is defined in the same way via the concordance number of the cycle.

If edges $e_i$ and $e_j$ in $P$ coincide and $i<j$, the path $P$ can be decomposed into the
path $P'=(e_1,\dots,e_{i-1},e_{i}=e_j,e_{j+1},\dots,e_{r})$
and the cycle $C^0=(e_{i},e_{i+1},\dots,e_{j-1})$. Clearly, $c(C_P)=c(C_{P'})+c(C^0)$, and hence 
\begin{equation}\label{offcycle}
w_P=-w_{P'}w_{C^0}.
\end{equation}

\begin{example} {\rm Consider the path $P=(e_1,e_2,e_3,e_{11},e_7,e_{10},e_3,e_{11},e_8)$ in 
Figure~\ref{fig:pergraph}. Choose $l$ as shown in the Figure; $l$ is neither collinear with the edges of $P$ nor with the relevant edges of the convex hull of boundary vertices (shown by dotted lines). Clearly,
$c_l(e',e'')=0$ for all pairs of consecutive edges of $C_P$ except for the pairs $(e_{10},e_3)$ and $(e_8,\bar e)$, where $\bar e$ is the additional edge joining $b_3$ and $b_4$. So, $c(C_P)=0$, and hence
$w_P=-w_1w_2w_3^2w_7w_8w_{10}w_{11}^2$. The same result can be obtained by decomposing $P$ into the path $P'=(e_1,e_2,e_3,e_{11},e_8)$ and the cycle $C^0=(e_3,e_{11},e_7,e_{10})$.}
\end{example}

\begin{remark}\label{dualdir}{\rm
Instead of closed polygonal curves $C_P$, one can use curves $C^*_P$ obtained by adding to $P$ the path between $b_j$ and $b_i$ that goes clockwise along the boundary of the convex hull of all the boundary vertices of $N$. It is a simple exercise to prove that the concordance numbers of $C_P$ and $C^*_P$ coincide. Therefore, the weight of a path can be defined also as 
\[
w_P=(-1)^{c(C^*_P)-1}\prod_{e\in P} w_e.
\]
}
\end{remark}

\subsection{Boundary measurements}\label{BM} Given a
perfect planar network as above, a sour\-ce $b_i$, $i\in I$, and a sink $b_j$, $j\in J$, we define the {\it boundary measurement\/} $M(i,j)$ as the
sum of the weights of all paths starting at $b_i$ and ending at $b_j$. Clearly, the boundary measurement thus defined is a formal infinite series
in variables $w_e$, $e\in E$. However, this series possesses certain nice propetsies.

Recall that a formal power series $g\in \Z[[w_e, e\in E]]$ is called a rational function if there exist polynomials $p,q\in \Z[w_e, e\in E]$ such that $p=qg$ in $\Z[[w_e, e\in E]]$. In this case we write $g=p/q$. For example, $1-z+z^2-z^3+\dots=(1+z)^{-1}$ in $\Z[[z]]$. Besides, we say that $g$ admits a {\it subtraction-free\/} rational expression if
it can be written as a ratio of two polynomials with nonnegative coefficients.
For example, $x^2-xy +y^2$ admits a subtraction-free rational expression since it can be written as $(x^3+y^3)/(x+y)$.

The following result was proved in \cite{Postnikov}[Lemma 4.3] and further generalized in \cite{Talaska}, but we will present an alternative proof
to illustrate the method that will be used in other proofs below.

\begin{proposition}\label{sfree}
Let $N$ be a perfect planar network in a disk, then each boundary measurement in $N$ is a rational function in the weights $w_e$ admitting a subtraction-free rational expression. 
\end{proposition}

\begin{figure}[ht]
\begin{center}
\includegraphics[height=8cm]{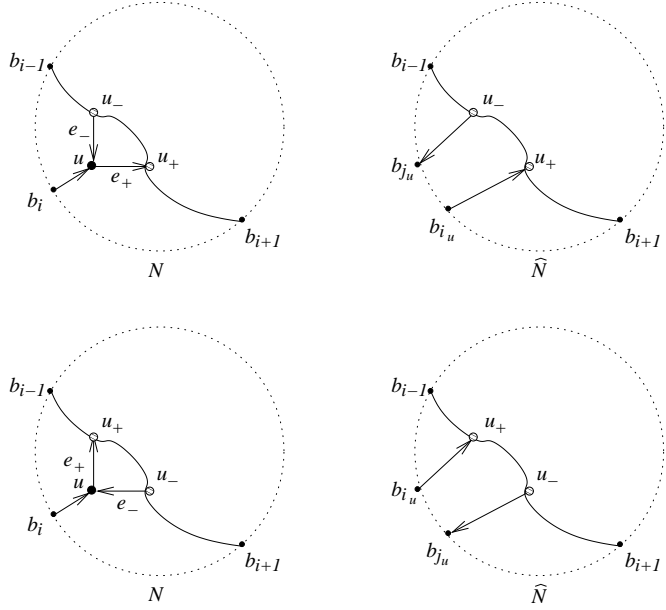}
\caption{Splitting a black vertex: cases $i-1\prec i_u\prec j_u\prec i+1$ 
(upper part) and $i-1\prec j_u\prec i_u\prec i+1$ (lower part)}
\label{fig:blackfork}
\end{center}
\end{figure}

\begin{proof}
We prove the claim by induction on the number of internal vertices. The base of
induction is the case when there are no internal vertices at all, and hence each
edge connects a source and a sink; in this case the statement of the proposition holds trivially.

Assume that $N$ has $r$ internal vertices. Consider a specific boundary
measurement $M(i,j)$.
The claim concerning $M(i,j)$ is trivial if
$b_j$ is the neighbor of $b_i$. In the remaining cases $b_i$
is connected by an edge $e_0$ to its only neighbor in $G$, which is either
white or black.

Assume first that the neighbor of $b_i$ is a white vertex $v$. Create a new network $\wN$ by
deleting $b_i$ and the edge $(b_i,v)$ from $G$, splitting $v$ into
$2$ sources $b_{i_v'}, b_{i_v''}$ (so that $i-1\prec i_v'\prec i_v''\prec i+1$ in the counterclockwise order) and replacing
the edges $e_1=(v,v')$ and $e_2=(v,v'')$ by $(b_{i_v'},v')$ and
$(b_{i_v''},v'')$, respectively, 
both of weight~1 (see Figure~\ref{fig:whitefork}).
Clearly, to any path $P$ from $b_i$ to $b_j$ corresponds either a path $P'$ from $b_{i'_v}$ to $b_j$ or a path $P''$ from $b_{i''_v}$ to $b_j$. Moreover, $c(C_P)=c(C_{P'})=c(C_{P''})$. Therefore
$$
M(i,j)=w_{e_0} (w_{e_1}\wM({i'_v},j)+w_{e_2}\wM({i''_v},j)),
$$
where $\wM$ means that the measurement is taken in $\wN$. Observe that the number
of internal vertices in $\wN$ is $r-1$, hence the claim follows from the above
relation by induction.

\begin{figure}[ht]
\begin{center}
\includegraphics[height=4.5cm]{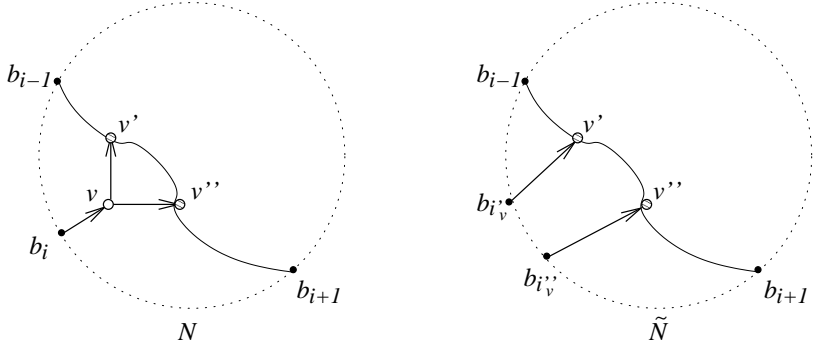}
\caption{Splitting a white vertex}
\label{fig:whitefork}
\end{center}
\end{figure}

Assume now that the neighbor of $b_i$ is a black vertex $u$. Denote by $u_+$ the unique vertex in
$G$ such that $(u,u_+)\in E$, and by $u_-$ the neighbor of $u$
distinct from $u_+$ and $b_i$. Create a new network $\hN$ by
deleting $b_i$ and the edge $(b_i,u)$ from $G$, splitting $u$ into
one new source $b_{i_u}$ and one new sink $b_{j_u}$ (so that either
$i-1\prec i_u\prec j_u\prec i+1$ or $i-1\prec j_u\prec
i_u\prec i+1$ in the counterclockwise order) and replacing the
edges $e_+=(u,u_+)$ and $e_-=(u_-,u)$ by new edges $\hat e_+$ and $\hat e_-$ in the same way as in the
previous case, see Figure~\ref{fig:blackfork}.

Let us classify the paths from $b_i$ to $b_j$ according to the number of times they traverse $e_+$. Similarly to the previous case, the total weight of the paths traversing $e_+$ only once is
given by $w_{e_0}w_{e_+}\hM({i_v},j)$, where $\hM$ means that the measurement is taken in $\hN$. Any path $P$ traversing $e_+$ exactly twice
can be represented as $P=(e_0,e_+,P_1,e_-,e_+,P_2)$ for some path $\hP_1=(\hat e_+,P_1,\hat e_-)$ from $i_u$ to $j_u$ in $\hN$ and a path $\hP_2=(\hat e_+,P_2)$ from $i_u$ to $b_j$ in $\hN$.
Clearly, $C_1=(e_+,P_1,e_-)$ is a cycle in $N$, and $c(C_1)=c(C_{\hP^1})$. Besides, $c(C_{P'})=c(C_{\hP_2})$ for $P'=(e_0,e_+,P_2)$. Therefore, by~(\ref{offcycle}), $w_P=-w_{C_1}w_{P'}$. Taking into account that
$w_{C_1}=w_{e_-}w_{e_+}w_{\hP_1}$ and $w_{P'}=w_{e_0}w_{e_+}w_{\hP_2}$,
we see that 
the total contribution of all paths traversing $e_+$ exactly twice to $M(i,j)$ equals
$$
-w_{e_0}w_{e_+}\hM({i_u},{j})w_{e_-}w_{e_+}\hM({i_u},j_u).
$$
In general, the total contribution of all paths traversing $e_+$ exactly $s+1$ times to $M(i,j)$ equals
$$
(-1)^{s}w_{e_0}w_{e_+}\hM({i_u},{j})(w_{e_-}w_{e_+}\hM({i_u},j_u))^s.
$$
Therefore, we find
\begin{equation}\label{sfree1}
M(i,j)=\frac{w_{e_0}w_{e_+}\hM({i_u},{j})}{1+w_{e_-}w_{e_+}\hM({i_u},j_u)}
\end{equation}
and the claim follows by induction, since the number of internal vertices in
$\hN$ is $r-1$.
\end{proof}

Boundary measurements can be organized into a $k\times m$ {\it
boundary measurement matrix\/} $M_N$ 
in the following way. Let $I=\{i_1<i_2<\dots<i_k\}$ and $J=\{j_1<j_2<\dots<j_m\}$.
We define $M_N=(M_{pq})$, $p\in [1,k]$, $q\in [1,m]$, where $M_{pq}=M({i_p},{j_q})$. 
Each network $N$ defines a rational map
$\CC_N\to \Mat_{k,m}$ given by $M_N$ and called the {\it boundary
measurement map\/} corresponding to $N$.  Here and below $ \Mat_{k,m}$ denotes the space
of $k\times m$ matrices.

\begin{example}\label{ex:bmmex}{\rm
Consider the network shown in Figure~\ref{fig:pergraph}. 
The corresponding boundary measurement matrix is a
$2\times 2$ matrix given by
\[
\begin{pmatrix}
\dfrac{w_3w_4w_5w_6w_{10}}{1+w_3w_7w_{10}w_{11}} &
\dfrac{w_3w_5w_6w_8w_{11}}{1+w_3w_7w_{10}w_{11}}\\
\\
\dfrac{w_1w_3w_4(w_2+w_6w_9w_{10})}{1+w_3w_7w_{10}w_{11}} &
\dfrac{w_1w_3w_8w_{11}(w_2+w_6w_9w_{10})}{1+w_3w_7w_{10}w_{11}}
\end{pmatrix}.
\]
}
\end{example}

\section{Poisson structures on the space of
edge weights and induced Poisson structures on
$\Mat_{k,m}$}\label{PSmat}

\subsection{Network concatenation and the standard Poisson-Lie structure}\label{concat}
A natural operation on networks is their {\it concatenation}, which consists, roughly
speaking, in gluing some sinks/sources of one network to  some of the
sources/sinks of the other.  We expect any Poisson structure associated
with networks to behave naturally under concatenation.
To obtain  from two planar networks a new one by concatenation, one needs 
to select a segment from the boundary of each disk
and identify these segments via a homeomorphism 
in such a way that
every sink (resp. source) contained in the selected  segment of the first
network is glued to a source (resp. sink)  of the second network. We can then erase the common piece of the boundary
along which the gluing was performed and identify every pair of glued edges
in the resulting network with a single edge of the same orientation and with the weight equal to the product
of two weights assigned to the two edges that were glued.

As an illustration, let us review a particular but important case, in which sources and sinks of the network do not interlace. In this case, it is more
convenient to view the network as located in a square rather than in a disk, with all sources located on the left side
and sinks on the right side of the square. It will be also handy to label
sources (resp. sinks) $1$  to $k$ (resp.~$1$  to $m$) going from the bottom to the top. This results in a different way
of recording boundary measurements into a matrix. Namely, if $M$ is the
boundary measurements matrix we defined earlier, then now we associate with the network the matrix $A = M W_0$, where $W_0=(\delta_{i,m+1-j})_{i,j=1}^m$ is the matrix of the longest permutation $w_0$.

\begin{figure}[ht]
\begin{center}
\includegraphics[height=3.0cm]{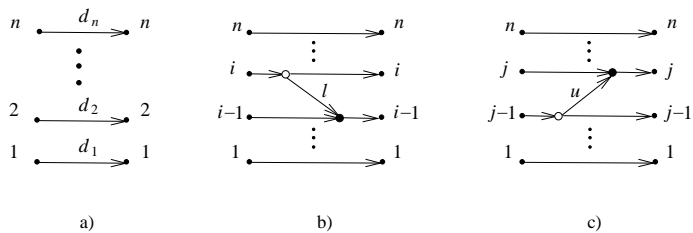}
\caption{Three networks used in matrix factorization}
\label{fig:factor}
\end{center}
\end{figure}

We can now concatenate
two networks of this kind, one with  $k$  sources and $m$ sinks and another
with $m$  sources and
$l$ sinks, by gluing the sinks of the former to the sources of the latter.
If $A_1$, $A_2$ are  $k\times m$ and $m\times l$ matrices associated with the two networks, then the matrix associated with their concatenation
is $A_1 A_2$. Note that this ``visualization'' of the matrix  multiplication
is particularly relevant when one deals with factorization of matrices into
products of elementary bidiagonal matrices.
Indeed, a
$n\times n$ diagonal  matrix $\diag(d_1,\ldots,d_n)$ and elementary
bidiagonal matrices $E^-_i(l) :=\one + l e_{i, i-1}$ and  $E^+_j(u) :=\one +
u e_{j-1, j}$ correspond to planar networks shown in Figure~\ref{fig:factor} a, b and c, respectively; all weights not shown explicitly are equal to~1.

The construction of Poisson structures we are about to present is motivated by
the way
in which the network representation of the bidiagonal factorization reflects
the {\it standard
Poisson-Lie structure}  on $SL_n$.
We need to recall some facts about Poisson-Lie groups
(see, e.g.\cite{r-sts}).

Let $\G$ be a Lie group equipped with a Poisson bracket $\Poi$.
$\G$ is called a {\em Poisson-Lie group}
if the multiplication map
$$
{\mathfrak m} : \G\times \G \ni (x,y) \mapsto x y \in \G
$$
is Poisson. Perhaps, the most important class of Poisson-Lie groups
is the one associated with classical R-matrices. 

Let $\g$ be a Lie algebra of
$\G$. Assume that $\g$ is equipped with a nondegenerate invariant
bilinear form $(\ ,\ )$. An element $R\in \End(\g)$ is a {\em classical R-matrix} if it is a skew-symmetric operator that satisfies the {\em modified classical
Yang-Baxter equation} (MCYBE)
\begin{equation}
[R(\xi), R(\eta)] - R \left ( [R(\xi), \eta]\ + \ [\xi, R(\eta)] \right ) = - [\xi,\eta]\ .
\label{MCYBE}
\end{equation}

Given a classical R-matrix $R$, $\G$ can be endowed with
a Poisson-Lie structure as follows.
Let $\nabla f, \nabla' f$ be the
right and the left gradients for a function $f\in C^\infty(\G)$:
\begin{equation}
( \nabla f ( x) , \xi ) = \frac{d} {dt} f( \exp{(t \xi)} x)\vert_{t=0},\qquad
 (\nabla' f ( x) , \xi ) = \frac{d} {dt} f( x\exp{(t \xi)})\vert_{t=0}\ .
\label{rightleftgrad}
\end{equation}
Then the  bracket given by 
\begin{equation}\nonumber
\{ f_1, f_2\} = \frac{1}{2} ( R(\nabla' f_1), \nabla' f_2 ) -
\frac{1}{2} ( R(\nabla f_1), \nabla f_2 )\  \label{sklyabra}
\end{equation}
is a Poisson-Lie bracket on $\G$ called the {\em Sklyanin bracket}.

We are interested in the case 
$\G=SL_n$ and $ \g=sl_n$ equipped with the trace-form
$$
(\xi, \eta) = \Tr ( \xi \eta)\ .
$$

Then the right and left gradients  (\ref{rightleftgrad}) are
$$
\nabla f(x) = x\ \grad f (x)\ , \qquad \nabla' f(x) = \grad f (x)\
x\ ,
$$
where
$$
\grad f (x) = \left (  \frac{\partial f} {\partial x_{ji}} \right
)_{i,j=1}^n\  ,
$$
and the Sklyanin bracket becomes
\begin{multline}
\{ f_1, f_2\}_{SL_n}(x) = \\
\frac{1}{2} ( R(\grad f_1(x)\ x), \grad
f_2(x)\ x ) -  \frac{1}{2} ( R(x\ \grad f_1(x)), x\ \grad f_2(x) )\ . \label{sklyaSLn}
\end{multline}

Every $\xi \in \mathfrak g$ can be uniquely decomposed as 
\begin{equation}\nonumber
\xi = \xi_- + \xi_0 + \xi_+\ ,
\label{decomposition_algebra}
\end{equation}
where $\xi_+$ and $\xi_-$ are strictly upper and lower triangular  and $\xi_0$
is diagonal. 
The simplest classical R-matrix on $sl_n$ is given by
\begin{equation}
R_0 (\xi) = \xi_+ - \xi_- = \left ( \s(j-i) \xi_{ij}\right )_{i,j=1}^n .
\label{standardR}
\end{equation}
 Substituting into (\ref{sklyaSLn}) $R=R_0$  we obtain , for matrix entries $x_{ij}, \
x_{i'j'}$, 
\begin{equation}\label{braijkl}
\{ x_{ij}, x_{i'j'}\}_{SL_n}=\frac{1}{2} \left ( \s(i'-i) +
\s(j'-j)\right ) x_{ij'} x_{i'j}\ .
\end{equation}

The bracket (\ref{sklyaSLn}) extends naturally to a Poisson bracket on the space
$\Mat_n$ of $n\times n$ matrices.

For example, the standard Poisson-Lie structure on
$$
SL_2 = \left \{ \left ( \begin{array} {cc} a & b\\ c & d
\end{array}\right )\ :\ ad-bc=1 \ \right \}
$$
is described by the relations
\begin{eqnarray*}
\nonumber
&\{ a, b\}_{SL_2} = \frac{1}{2} a b,\quad \{ a, c\}_{SL_2} = \frac{1}{2} a
c,
\quad \{ a, d\}_{SL_2} =  b c, \\
\nonumber
&\{ c, d\}_{SL_2} = \frac{1}{2} c d,\quad \{ b, d\}_{SL_2} = \frac{1}{2} b
d,
\quad \{ b, c\}_{SL_2} =  0\ ,
\end{eqnarray*}
which, when restricted to  upper and lower Borel subgroups  of $
SL_2$
\[
B_+=\left\{
\begin{pmatrix} d  & c \\ 0 & d^{-1}
\end{pmatrix}\right \}, \qquad
B_-=\left  \{ \begin{pmatrix} d  & 0 \\
c & d^{-1} \end{pmatrix}\right \}
\]
have an especially simple 
form
\[
\{ d, c\} = \frac{1}{2} dc.
\]
The latter Poisson brackets can be used to give an alternative characterization
of the standard Poisson-Lie structure on $SL_n$. Namely, define  the canonical embedding
$\rho_i: SL_2 \to SL_n $ ($i\in [1,
n-1]$) 
that maps $SL_2$ into $SL_n$ as a diagonal
$2\times 2$ block occupying rows and columns $i$ and $i+1$. 
Then the standard Poisson-Lie structure on
$SL_n$ is defined uniquely (up to a scalar multiple) by the requirement
that restrictions of $\rho_i$ to $B_\pm$ are Poisson.

\begin{figure}[ht]
\begin{center}
\includegraphics[height=3.0cm]{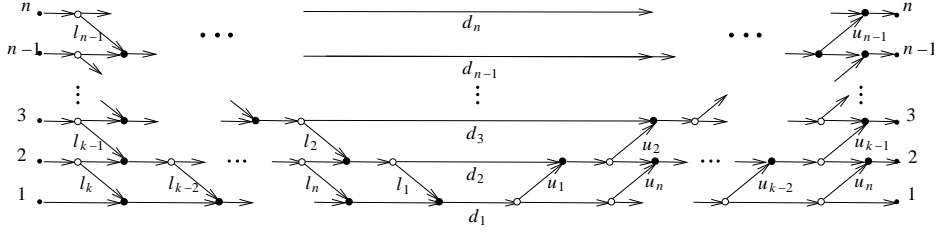}
\caption{Generic planar network}
\label{fig:genfactor}
\end{center}
\end{figure}

Note that the network that represents $\rho_i(B_-)$ looks like the second
network in the figure above
with the weights $d, d^{-1}$ and $c$ attached to edges $(i-1) \to (i-1)$, $i
\to i$ and $i \to (i-1)$ resp.,
while the network that represents $\rho_j(B_+)$ looks like the third network
in the figure above
with the weights $d, d^{-1}$ and $c$ attached to edges $(j-1) \to (j-1)$,
$j\to j$ and $(j-1) \to j$ resp.
Concatenation of several networks  $N_1, \cdots, N_r$, $r=n(n-1)$, of these two kinds,
with appropriately chosen order and with each diagram having its own pair of
nontrivial weights $c_i,d_i$, describes a generic
element of $SL_n$ (see, e.g. \cite{Fallat}). An example of such a network is given
by Figure~\ref{fig:genfactor}.


On the other hand, due to the Poisson-Lie property, the
Poisson structure on $SL_n$ is inherited from simple Poisson brackets
for parameters $c_i, d_i$, which can be described completely in terms of
networks: (i)  the bracket of any two parameters is equal to their product times a constant; (ii) this constant is equal to zero unless the corresponding  edges  have a
common source/sink; (iii) the constant
 equals $\pm \frac{1}{2}$ if
the corresponding edges follow one another around the source/sink in the
counterclockwise direction.
The corresponding Poisson-Lie structure on $GL_n$ is obtained by requiring
the determinant to be a Casimir function.

This example motivates conditions we impose below on a natural Poisson
structure associated
with a 3-valent planar directed network.

\subsection{Poisson structures on the space of
edge weights}\label{PSSC} 
\begin{figure}[ht]
\begin{center}
\includegraphics[height=1.5cm]{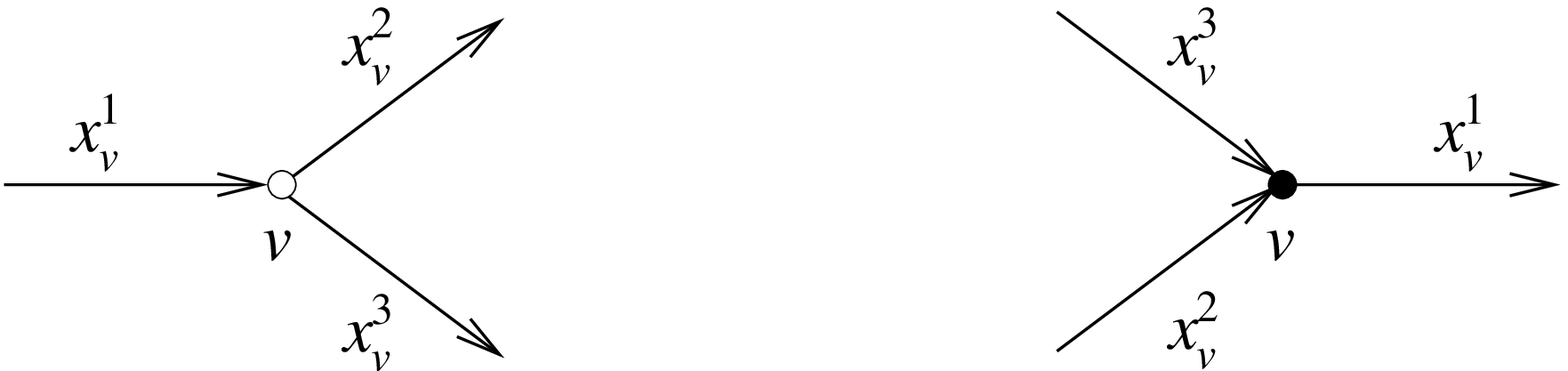}
\caption{Edge labelling for $\RS_v^3$} \label{fig:trefoils}
\end{center}
\end{figure}

 Let G be a directed planar graph in a disk as
described in Section~\ref{PandW}.
A pair $(v,e)$ is called a {\it flag\/} if $v$ is an endpoint of $e$.
To each internal vertex $v$ of $G$ we
assign a 3-dimensional space $\RS_v^3$ with coordinates $x_v^1,
x_v^2, x_v^3$. We equip each $\RS_v^3$ with a Poisson bracket
$\{\cdot,\cdot\}_v$. It is convenient to assume that the 
flags involving $v$ 
 are labelled by the coordinates, as shown in Figure~\ref{fig:trefoils}.

Besides, to each boundary vertex $b_j$ of $G$ we assign a
1-dimensional space $\RS_j$ with the coordinate $x_j^1$ (in accordance with the above convention, this coordinate labels the 
unique flag involving
$b_j$). Define $\RR$
to be the direct sum of all the above spaces; thus, the dimension
of $\RR$ equals twice the number of edges in $G$. Note that $\RR$ is
equipped with a Poisson bracket $\{\cdot,\cdot\}_\RR$, which is
defined as the direct sum of the brackets $\{\cdot,\cdot\}_v$; that is, $\{x,y\}_\RR=0$ whenever $x$ and $y$ are not defined on the same $\RS_v^3$. We
say that the bracket $\{\cdot,\cdot\}_\RR$ is {\it universal\/} if
each of $\{\cdot,\cdot\}_v$ depends only on the color of the
vertex $v$.

Define the weights $w_e$ by 
\begin{equation}\label{connect}
w_e=x^i_vx^j_u, 
\end{equation}
provided the 
flag $(v,e)$ is labelled by $x^i_v$ and the flag $(u,e)$ is labelled by $x^j_u$.
In other words, the weight of an edge is defined as the product of the weights of the two 
flags involving 
this edge.
Therefore, in this case the space of edge weights $\CC_N$ coincides with the entire $\RS^{|E|}$, and the weights define a {\it weight map\/} $w\: \RS^d\to\RS^{|E|}$. 
We require the pushforward of 
$\{\cdot,\cdot\}_\RR$ to $\RS^{|E|}$ by the weight map to be a well defined 
Poisson bracket; this can be regarded as an analog
of the Poisson--Lie property for groups.

\begin{proposition}\label{6param}
Universal Poisson brackets $\{\cdot,\cdot\}_\RR$ such that the
weight map $w$ is Poisson form a $6$-parametric family defined
by relations
\begin{equation}\label{6parw}
\{x^i_v,x^j_v\}_v=\alpha_{ij}x^i_vx^j_v, \quad i,j\in [1,3], i\ne
j,
\end{equation}
at each white vertex $v$ and
\begin{equation}\label{6parb}
\{x^i_v,x^j_v\}_v=\beta_{ij}x^i_vx^j_v, \quad i,j\in [1,3], i\ne
j,
\end{equation}
at each black vertex $v$.
\end{proposition}

\begin{proof} Indeed, let $v$ be a white vertex, and let $e=(v,u)$ and
$\bar e=(v,\bar u)$ be the two outcoming edges. By definition,
there exist $i,j,k,l\in [1,3]$, $i\ne j$, such that
$w_e=x_v^ix_u^k$, $w_{\bar e}=x_v^jx_{\bar u}^l$. Therefore,
$$
\{w_e,w_{\bar e}\}_N=\{x_v^ix_u^k,x_v^jx_{\bar u}^l\}_\RR=x_u^kx_{\bar
u}^l\{x_v^i,x_v^j\}_v,
$$
where $\{\cdot,\cdot\}_N$ stands for the pushforward of
$\{\cdot,\cdot\}_\RR$. Recall that the Poisson bracket in $\RS_v^3$
depends only on $x_v^1$, $x_v^2$ and $x_v^3$. Hence the only
possibility for the right hand side of the above relation to be a
function of $w_e$ and $w_{\bar e}$ occurs when
$\{x^i_v,x^j_v\}_v=\alpha_{ij}x^i_vx^j_v$, as required.

Black vertices are treated in the same way.
\end{proof}

Let $v$ be a white vertex. A {\it local gauge
transformation\/} at $v$ is a transformation $\RS_v^3\to\RS_v^3$
defined by $(x_v^1, x_v^2, x_v^3)\mapsto (\bar x_v^1=x_v^1t_v, \bar x_v^2=x_v^2t_v^{-1}, 
\bar x_v^3=x_v^3t_v^{-1})$, where $t_v$ is a Laurent monomial in $x_v^1$,
$x_v^2$, $x_v^3$.
A local gauge transformation at a black vertex is defined by the
same formulas, with $t_v$ replaced by $t_v^{-1}$.

A {\it global gauge transformation\/} $t\:\RR\to \RR$ is defined by
applying a local gauge transformation $t_v$ at each vertex $v$.
The composition map $w\circ t$ defines a network $tN$; the graph
of $tN$ coincides with the graph of $N$, and the weight $w^t_e$ of
an edge $e=(u,v)$ is given by $w^t_e=t_vw_et_u^{-1}$. Therefore,
the weights of the same path in $N$ and $tN$ coincide. It follows
immediately that
\begin{equation}\label{gaugcom}
M_{tN}\circ w\circ t=M_N\circ w,
\end{equation}
provided both sides of the equality are well defined.

\subsection{Induced Poisson structures on
$\Mat_{k,m}$}\label{PSonmat}
Our next goal is to look at Poisson properties of the boundary
measurement map. Fix an arbitrary partition $I\cup J=[1,n]$, $I\cap J=\varnothing$,
and let $k=|I|$, $m=n-k=|J|$.
Let $\Net_{I,J}$ stand for the set of all perfect planar
networks in a disk with the sources $b_i$, $i\in I$, sinks $b_j$, $j\in J$, and
edge weights $w_e$ defined by~(\ref{connect}).
We assume that the space of edge weights $\CC_N=\RS^{|E|}$ is
equipped with the Poisson bracket $\{\cdot,\cdot\}_N$ obtained as the pushforward of the 6-parametric family $\{\cdot,\cdot\}_\RR$ described in
Proposition~\ref{6param}.

\begin{theorem}\label{PSM} 
There exists a $2$-parametric family of Poisson brackets on
$\Mat_{k,m}$ with the following property: for any choice of
parameters $\alpha_{ij}$, $\beta_{ij}$
in~\eqref{6parw},~\eqref{6parb} this family contains a unique
Poisson bracket on $\Mat_{k,m}$ such that for any network
$N\in\Net_{I,J}$ the map $M_N\:\RS^{|E|}\to \Mat_{k,m}$ is Poisson.
\end{theorem}

\begin{proof} 
Relation~(\ref{gaugcom}) suggests that one may use global gauge transformations in order to decrease the number of
parameters in the universal 6-parametric family described in
Proposition~\ref{6param}.
Indeed, for any white vertex $v$ we consider a local gauge
transformation $(x_v^1,x_v^2,x_v^3)\mapsto (\bar x_v^1,\bar x_v^2,
\bar x_v^3)$ with $t_v=1/x_v^1$. Evidently,
\begin{equation}\label{2parw}
\{\bar x^2_v,\bar x^3_v\}_v=\alpha \bar x^2_v\bar x^3_v, \quad
\{\bar x^1_v,\bar x^2_v\}_v=\{\bar x^1_v,\bar x^3_v\}_v=0
\end{equation}
with 
\begin{equation}\label{condw}
\alpha=\alpha_{23}+\alpha_{13}-\alpha_{12}.
\end{equation}

 Similarly, for
each black vertex $v$ we consider a local gauge transformation
$(x_v^1,x_v^2,x_v^3)\mapsto (\bar x_v^1,\bar x_v^2, \bar x_v^3)$
with $t_v=x_v^1$. Evidently,
\begin{equation}\label{2parb}
\{\bar x^2_v,\bar x^3_v\}_v=\beta \bar x^2_v\bar x^3_v, \quad
\{\bar x^1_v,\bar x^2_v\}_v=\{\bar x^1_v,\bar x^3_v\}_v=0
\end{equation}
with 
\begin{equation}\label{condb}
\beta=\beta_{23}+\beta_{13}-\beta_{12}.
\end{equation}

From now on we consider the 2-parametric
family~(\ref{2parw}),~(\ref{2parb}) instead of the 6-parametric
family~(\ref{6parw}),~(\ref{6parb}). 

To define a Poisson bracket on $\Mat_{k,m}$, it suffices to
calculate the bracket for any pair of matrix entries and to extend it further via bilinearity and the Leibnitz identity. To do this we will need the following two auxiliary functions: for any $i,i'\in I$, $j,j' \in J$ define
\begin{equation}\label{spar}
s_=(i,j,i',j')=\begin{cases} 
1\qquad\text{if $i\prec i'\prec j'\prec j$}, \\
-1\qquad\text{if $i'\prec i\prec j\prec j'$},\\
\frac12\qquad\text{if $i=i'\prec j'\prec j$ or $i\prec i'\prec j'=j$}, \\ 
-\frac12\qquad\text{if $ i'=i\prec j\prec j'$ or $i'\prec i\prec j=j'$},\\
0 \qquad \text{otherwise},
\end{cases}
\end{equation} 
and
\begin{equation}\label{scross}
s_\times(i,j,i',j')=\begin{cases} 
1\qquad\text{if $i'\prec i\prec j'\prec j$}, \\
-1\qquad\text{if $i\prec i'\prec j\prec j'$},\\
\frac12\qquad\text{if  $i'=i\prec j'\prec j$ or $i'\prec i\prec j'=j$}, \\ 
-\frac12\qquad\text{if $i=i'\prec j\prec j'$ or $i\prec i'\prec j=j'$},\\
0 \qquad \text{otherwise}.
\end{cases}
\end{equation}  
Note that both $s_=$ and $s_\times$ are skew-symmetric:
\[
s_=(i,j,i',j')+s_=(i',j',i,j)=s_\times(i,j,i',j')+s_\times(i',j',i,j)=0
\]
for any $i,i'\in I$, $j,j'\in J$. Quadruples $(i,j,i',j')$ such that at least one of the values $s_=(i,j,i',j')$ and $s_\times(i,j,i',j')$ is distinct from zero are shown in Figure~\ref{fig:twofunctions}. For a better visualisation, pairs $i,j$ and $i',j'$ are joined by a directed edge; these edges should not be mistaken for edges of $N$.

\begin{figure}[ht]
\begin{center}
\includegraphics[height=5.5cm]{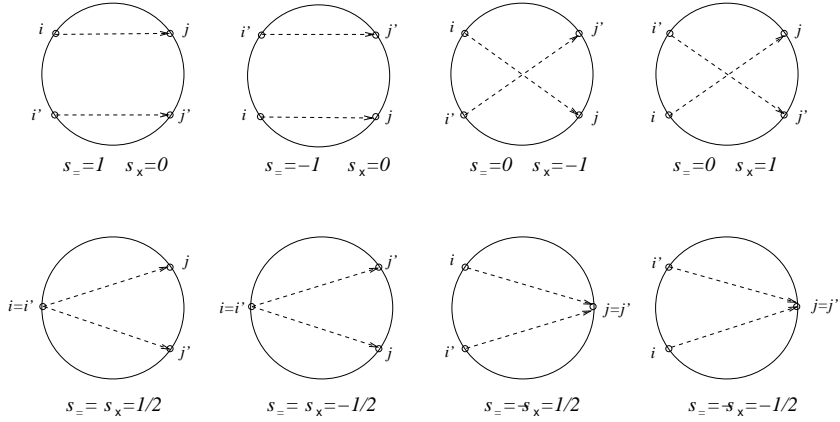}
\caption{Nontrivial values of $s_=$ and $s_\times$} \label{fig:twofunctions}
\end{center}
\end{figure}

Theorem~\ref{PSM} is proved by presenting an explicit formula for the bracket on $\Mat_{k,m}$.
\end{proof}
 
\begin{theorem}\label{PSME}
The $2$-parametric family of Poisson brackets on $\Mat_{k,m}$ satisfying 
the conditions of Theorem~\ref{PSM} is given by
\begin{equation}\label{psme}
\{M_{pq},M_{\bar p\bar q}\}_{I,J}=
(\alpha-\beta)s_=(i_p,j_q,i_{\bar p},j_{\bar q})M_{p\bar q}M_{\bar p q}+
(\alpha+\beta)s_\times(i_p,j_q,i_{\bar p},j_{\bar q})M_{pq}M_{\bar p\bar q},
\end{equation}
where $\alpha$ and $\beta$ satisfy~\eqref{condw},~\eqref{condb},
$p,\bar p\in [1,k]$, $q,\bar q\in [1,m]$.
\end{theorem}

\begin{proof}
First of all, let us check that that relations~(\ref{psme}) indeed define a Poisson bracket on $\Mat_{k,m}$. Since bilinearity and the Leibnitz identity are built-in in the definition, and skew symmetry follows immediately from~(\ref{spar}) and~(\ref{scross}), it remains to check the Jacobi identity. 

\begin{lemma}\label{jacobi}
The bracket $\{\cdot,\cdot\}_{I,J}$ satisfies the Jacobi identity.
\end{lemma}

\begin{proof}
The claim can be verified easily when at least one of the following five conditions holds true: $i=i'=i''$; $j=j'=j''$; $i=i'$ and $j=j'$; $i=i''$ and $j=j''$; $i'=i''$ and $j'=j''$. In what follows we assume that none of these conditions holds.

A simple computation shows that under this assumption, the Jacobi identity for 
$\{\cdot,\cdot\}_{I,J}$ is implied by the following three identities for the functions $s_=$ and $s_\times$: for any $i,i',i''\in I$, $j,j',j''\in J$,
\begin{align}
s_\times&(i',j',i'',j'')s_\times(i,j,i'',j'')+
s_\times(i'',j'',i,j)s_\times(i',j',i,j)\notag\\
&+s_\times(i,j,i',j')s_\times(i'',j'',i',j')+
s_\times(i',j',i'',j'')s_\times(i,j,i',j')\notag\\
&+s_\times(i'',j'',i,j)s_\times(i',j',i'',j'')+
s_\times(i,j,i',j')s_\times(i'',j'',i,j)=0, \notag\\
s_=&(i',j',i'',j'')s_=(i,j,i'',j')+s_=(i'',j'',i,j)s_=(i',j',i,j'')\label{threeid}\\
&+s_=(i,j,i',j')s_=(i'',j'',i',j)=0,\notag\\
s_=&(i',j',i'',j'')
(s_\times(i,j,i'',j')+s_\times(i,j,i',j'')-s_\times(i,j,i'',j'')\notag\\
&-s_\times(i,j,i,'j'))=0.\notag
\end{align}

The first identity in~(\ref{threeid}) is evident, since by~(\ref{scross}), the first term is canceled by the fifth one, the second term is canceled by the sixth one, and the third term is canceled by the fourth one. 

To prove the second identity, assume to the contrary that there exist $i,i',i''\in I$, $j,j',j''\in J$ such that the left hand side does not vanish. Consequently, at least one of the three terms in the left hand side does not vanish; without loss of generality we may assume that it is the first term. 

Since $s_=(i',j',i'',j'')\ne0$, we get either 
\begin{equation}\label{order1}
i'\preceq i''\prec j'' \preceq j', 
\end{equation}
or  
\begin{equation}\label{order2}
i'\succeq i''\succ j'' \succeq j'.
\end{equation}

Assume that~(\ref{order1}) holds, than we have five possibilities for $j$:

(i) $j''\prec j\prec j'$; 

(ii) $j=j''$;

(iii) $i''\prec j \prec j''$;

(iv) $j'\prec j \prec i''$;

(v) $j=j'$.

In cases (i)-(iii) condition  $s_=(i,j,i'',j')\ne 0$ implies $i''\preceq i \prec j$. Therefore, in case (i) we get
\begin{align*}
s_=(i',j',i'',j'')=-&s_=(i'',j'',i',j), \qquad
s_=(i,j,i'',j')=s_=(i,j,i',j'), \\ 
&s_=(i'',j'',i,j)=0
\end{align*}
provided $i''\ne i$, and  
\begin{align*}
s_=(i',j',i'',j'')=-&s_=(i'',j'',i',j)=
s_=(i',j',i,j'')=-s_=(i,j,i',j')=1, \\
&s_=(i,j,i'',j')=s_=(i'',j'',i,j)=-1/2
\end{align*}
provided $i''=i\ne i'$. In both situations the second identity in~(\ref{threeid}) follows immediately.

In case (ii) we get
\[
s_=(i,j,i'',j')=s_=(i,j,i',j')=-s_=(i',j',i,j'')=-1,\qquad
s_=(i'',j'',i,j)=1/2,
\]
and 
\[
s_=(i',j',i'',j'')=\begin{cases} 1\quad&\text{if $i'\ne i''$},\\
                                 1/2\quad &\text{if $i'=i''$},
                   \end{cases}\qquad
s_=(i'',j'',i',j)  =\begin{cases} 1/2\quad&\text{if $i'\ne i''$},\\
                                 0\quad &\text{if $i'=i''$}.
                   \end{cases}                               
\]
In both situations the second identity in~(\ref{threeid}) follows immediately.

In case (iii) we get
\begin{align*}
s_=(i',j',i'',j'')=&s_=(i',j',i,j''), \qquad
s_=(i,j,i'',j')=-s_=(i'',j'',i,j), \\
&s_=(i'',j'',i',j)=0
\end{align*}
provided $i'\ne i''$, and  
\begin{align*}
s_=(i'',j'',i,j)=-&s_=(i,j,i'',j')=
s_=(i',j',i,j'')=-s_=(i,j,i',j')=1, \\ 
&s_=(i',j',i'',j'')=s_=(i'',j'',i',j)=1/2
\end{align*}
provided $i'=i''\ne i$. In both situations the second identity in~(\ref{threeid}) follows immediately.

Case (iv), in its turn, falls into three cases depending on the location of $i$; these three cases are parallel to the cases (i)-(iii) above and are treated in the same way. 

Finally, in case (v) we have to distinguish two subcases: $i''\prec i \prec i'$ and $i'\prec i\preceq i''$. In the first subcase we have
\begin{align*}
s_=(i',j',i'',j'')=-&s_=(i'',j'',i',j), \qquad
s_=(i,j,i'',j')=s_=(i,j,i',j'), \\ 
&s_=(i'',j'',i,j)s_=(i',j',i,j'')=0,
\end{align*}
while in the second subcase,
\[
s_=(i',j',i'',j'')=s_=(i',j',i,j'')=-s_=(i'',j'',i',j)=-1,\quad
s_=(i,j,i',j')=-1/2,
\]
and 
\[
s_=(i'',j'',i,j)=\begin{cases} -1\quad&\text{if $i\ne i''$},\\
                                -1/2\quad &\text{if $i=i''$},
                   \end{cases}\qquad
s_=(i,j,i'',j')  =\begin{cases} 1/2\quad&\text{if $i\ne i''$},\\
                                 0\quad &\text{if $i=i''$}.
                   \end{cases}                               
\]
In both situations the second identity in~(\ref{threeid}) follows immediately.

If the points are ordered counterclockwise as prescribed by~(\ref{order2}), the proof is very similar, with $\prec$ and $\preceq$ replaced by $\succ$ and $\succeq$, correspondingly.

To prove the third identity in~(\ref{threeid}), assume to the contrary that there exist $i,i',i''\in I$, $j,j',j''\in J$ such that the left hand side does not vanish. Consequently, $s_=(i',j',i'',j'')\ne 0$, and hence once again one of~(\ref{order1}) and~(\ref{order2}) holds.

\begin{figure}[ht]
\begin{center}
\includegraphics[height=2.5cm]{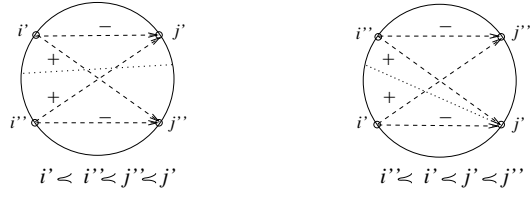}
\caption{To the proof of Lemma~\ref{jacobi}} \label{fig:btfly}
\end{center}
\end{figure}

Denote by $S$ the sum in the left hand side of the identity.
If $i'=i''$ or $j'=j''$ then $S$ vanishes due to the skew-symmetry of $s_\times$. The remaining case amounts to checking all possible ways to insert a chord $ij$ in the configurations presented in Figure~\ref{fig:btfly}. The check itself in each case is trivial. For example, if the inserted chord is as shown by the dotted line in the left part of Figure~\ref{fig:btfly}, then the last two terms of $S$ vanish, while the first two terms have opposite signs and absolute value~1. If the inserted chord is as shown by a dotted line in the right part of Figure~\ref{fig:btfly}, then one of the last two terms in $S$ vanishes and one of the first two terms has absolute value~1. The two remaining terms have absolute value $1/2$; besides, they have the same sign, which is opposite to the sign of the term with absolute value~1. Other cases are similar and left to the reader.
\end{proof} 

To complete the proof of Theorem~\ref{PSME}, it remains to check that
\begin{equation}\label{twops}
\{M_{pq},M_{\bar p\bar q}\}_N=\{M_{pq},M_{\bar p\bar q}\}_{I,J}
\end{equation}
for any pair of matrix entries $M_{pq}$ and $M_{\bar p\bar q}$. 
The proof is similar to the proof of Proposition~\ref{sfree} and relies on the induction on the number of inner vertices in $N\in \Net_{I,J}$.

Assume first that $N$ does not have inner vertices, and hence each edge of $N$
connects two boundary vertices. It is easy to see that in this case the Poisson bracket computed in $\RS^{|E|}$ vanishes identically.
Let us prove that the bracket given by~(\ref{psme}) vanishes as well. 

We start with the case when both $(b_{i_p},b_{j_q})$ and $(b_{i_{\bar p}},b_{j_{\bar q}})$ are edges in $N$. Then $p\ne \bar p$ and $q\ne\bar q$, since there is only one edge incident to each boundary vertex. Therefore
$M_{p\bar q}=M_{\bar pq}=0$, and the first term in the right hand side of~(\ref{psme}) vanishes. Besides, since $N$ is planar, 
$s_\times(i_p,j_q,i_{\bar p},j_{\bar q})=0$, and the second term vanishes as well. 

Next, let $(b_{i_p},b_{j_q})$ be an edge of $N$, and $(b_{i_{\bar p}},b_{j_{\bar q}})$ be a non-edge. Then $M_{pq}=M_{p\bar q}=0$, and hence both terms in the right hand side of~(\ref{psme}) vanish.

Finally, let both $(b_{i_p},b_{j_q})$ and $(b_{i_{\bar p}},b_{j_{\bar q}})$ be non-edges. Then $M_{pq}=M_{\bar p\bar q}=0$, and the second term in the right hand side of~(\ref{psme}) vanishes. The first term can be distinct from zero only if both $(b_{i_p},b_{j_{\bar q}})$ and $(b_{i_{\bar p}},b_{j_q})$ are edges in $N$. Once again we use planarity of $N$ to see that in this case 
$s_=(i_p,j_q,i_{\bar p},j_{\bar q})=0$, and hence the right hand side of~(\ref{psme}) vanishes.

Now we may assume that $N$ has $r$ inner vertices, and that~(\ref{twops}) is true for all networks with at most $r-1$ inner vertices and any number of boundary vertices. 
Consider the unique neighbor of $b_{i_p}$ in $N$. If this neighbor is another boundary vertex then the same reasoning as above applies to show that 
$\{M_{pq},M_{\bar p\bar q}\}_{I,J}$ vanishes identically for any choice of $i_{\bar p}\in I$, $j_q,j_{\bar q}\in J$, which agrees with the behavior of $\{\cdot,\cdot\}_N$.

Assume that the only neighbor of $b_{i_p}$ is a white inner vertex $v$. 
Define
\[
\widetilde{\RR}=\left(\RR\oplus\RS_{i'_v}\oplus\RS_{i''_v}\right)\ominus\left(\RS_v^3\oplus\RS_{i_p}\right),
\]
which corresponds to a network $\wN$ obtained from $N$ by deleting vertices $b_{i_p}$, $v$ and the edge $(b_{i_p},v)$ from $G$ and adding two new sources $i'_v$ and $i''_v$ so that $i_p-1\prec i'_v\prec i''_v\prec i_p+1$, see Figure~\ref{fig:whitefork}. Let $1, x', x''$ be the coordinates in $\RS_v^3$ after the local gauge transformation at $v$, so that $\{x',x''\}_v=\alpha x'x''$. Then 
\begin{equation}\label{simplerecalc}
M({i_p},j)=x'\wM({i'_v},j)+x''\wM({i''_v},j),\qquad
M({i_{\bar p}},j)=\wM({i_{\bar p}},j)
\end{equation}
for any $j\in J$ and any $\bar p\ne p$. Since $\wN$ has $r-1$ inner vertices, 
$\wM(i,j)$ satisfy relations similar to~(\ref{twops}) with $N$ replaced by $\wN$ and $I$ replaced by $I\setminus i_p\cup i'_v\cup i''_v$.
Relations~(\ref{twops}) follow immediately from this fact and~(\ref{simplerecalc}), provided $\bar p\ne p$. In the latter case we have
\begin{align*}
\{M_{pq},M_{p\bar q}\}_N
&=\{x'\wM({i'_v},{j_q})+x''\wM({i''_v},{j_q}),
x'\wM({i'_v},{j_{\bar q}})+x''\wM({i''_v},{j_{\bar q}})\}_N\\
&=(x')^2\{\wM({i'_v},{j_q}),\wM({i'_v},{j_{\bar q}})\}_\wN
+x''x'\{\wM({i''_v},{j_q}),\wM({i'_v},{j_{\bar q}})\}_\wN\\
&+\{x'',x'\}_v\wM({i''_v},{j_q})\wM({i'_v},{j_{\bar q}})
+x'x''\{\wM({i'_v},{j_q}),\wM({i''_v},{j_{\bar q}})\}_\wN\\
&+\{x',x''\}_v\wM({i'_v},{j_q})\wM({i''_v},{j_{\bar q}})
+(x'')^2\{\wM({i''_v},{j_q}),\wM({i''_v},{j_{\bar q}})\}_\wN.
\end{align*}

The first term in the right hand side of the expression above equals
\[
\left((\alpha-\beta)s_=(i'_v,j_q,i'_v,j_{\bar q})
+(\alpha+\beta)s_\times(i'_v,j_q,i'_v,j_{\bar q})\right)
x'\wM({i'_v},{j_q})x''\wM({i'_v},{j_{\bar q}}).
\]
Since $s_=(i'_v,j_q,i'_v,j_{\bar q})=s_\times(i'_v,j_q,i'_v,j_{\bar q})=
\pm\frac12$ (the sign is negative if $i_p\prec j_q\prec j_{\bar q}$ and positive if $i_p\prec j_{\bar q}\prec j_q$), the first term equals 
$\pm\alpha x'\wM({i'_v},{j_q})x''\wM({i'_v},{j_{\bar q}})$.

Similarly, the second term equals
\begin{align*}
-&(\alpha-\beta)x'\wM({i'_v},{j_q})x''\wM({i''_v},{j_{\bar q}})
\qquad\text{if $i_p\prec j_q\prec j_{\bar q}$}\\
&(\alpha+\beta)x'\wM({i'_v},{j_{\bar q}})x''\wM({i''_v},{j_q})
\qquad\text{if $i_p\prec j_{\bar q}\prec j_q$},
\end{align*}
the fourth term equals
\begin{align*}
-&(\alpha+\beta)x'\wM({i'_v},{j_q})x''\wM({i''_v},{j_{\bar q}})
\qquad\text{if $i_p\prec j_q\prec j_{\bar q}$}\\
&(\alpha-\beta)x'\wM({i'_v},{j_{\bar q}})x''\wM({i''_v},{j_q})
\qquad\text{if $i_p\prec j_{\bar q}\prec j_q$},
\end{align*}
and the sixths term equals $\pm\alpha x''\wM({i''_v},{j_q})x''\wM({i''_v},{j_{\bar q}})$ with the same sign rule as for the first term. We thus see that
\begin{align*}
\{M_{pq},M_{p\bar q}\}_N&=
\pm\alpha (x'\wM({i'_v},{j_q})+x''\wM({i''_v},{j_q}))
(x'\wM({i'_v},{j_{\bar q}})+x''\wM({i''_v},{j_{\bar q}}))\\
&=\pm\alpha M_{pq}M_{p\bar q}=\{M_{pq},M_{p\bar q}\}_{I,J},
\end{align*}
since $s_=(i_p,j_q,i_p,j_{\bar q})=s_\times(i_p,j_q,i_p,j_{\bar q})=
\pm\frac12$.  

Assume now that the only neighbor of $b_{i_p}$ is a black inner vertex $u$. 
Define
\[
\widehat{\RR}=\left(\RR\oplus\RS_{i_u}\oplus\RS_{j_u}\right)\ominus\left(\RS_u^3\oplus\RS_{i_p}\right),
\]
which corresponds to a network $\hN$ obtained from $N$ by deleting vertices $b_{i_p}$, $u$ and the edge $(b_{i_p},u)$ from $G$ and adding a new source $i_u$ and a new sink $j_u$ so that either $i_p-1\prec i_u\prec j_u\prec i_p+1$, or $i_p-1\prec i_u\prec j_u\prec i_p+1$, 
see Figure~\ref{fig:blackfork}. Let $1, x', x''$ be the coordinates in $\RS_u^3$ after the local gauge transformation at $u$, so that $\{x',x''\}_u=\beta x'x''$. 

\begin{lemma}\label{recalc}
Boundary measurements in the networks $N$ and $\hN$ are related by
\begin{align*}
M({i_p},j)&=\frac{x'\hM({i_u},j)}{1+x''\hM({i_u},{j_u})},\\
M({i_{\bar p}},j)&=\hM({i_{\bar p}},j)\pm
\frac{x''\hM({i_{\bar p}},{j_u})\hM({i_u},j)}{1+x''\hM({i_u},{j_u})},\qquad
\bar p\ne p;
\end{align*}
in the second formula above, sign $+$ corresponds to the cases 
\[
i_p-1\prec j_u\prec i_u\prec i_p+1\preceq j \prec i_{\bar p}
\]
or
\[
i_{\bar p}\prec j\preceq i_p-1\prec i_u\prec j_u\prec i_p+1,
\] 
and sign $-$ corresponds to the cases 
\[
i_p-1\prec i_u\prec j_u\prec i_p+1\preceq j \prec i_{\bar p}
\]
or
\[
i_{\bar p}\prec j\preceq i_p-1\prec j_u\prec i_u\prec i_p+1.
\]
\end{lemma}

\begin{proof}
The first formula above was, in fact, already obtained in the proof of Proposition~\ref{sfree}; one has to take into account that after the local gauge 
transformation at $u$ we get $w_{e_+}=1$, $w_{e_-}=x''$ and $w_{e_0}=x'$.
 
To get the second formula, we apply the same reasoning as in the proof of Proposition~\ref{sfree}. The paths from $b_{i_{\bar p}}$ to $b_j$ in $N$ are classified according to the number of times they traverse the edge $e_+$. 
The total contribution of the paths not traversing $e_+$ at all 
to $M({i_{\bar p}},j)$ equals
$\hM({i_{\bar p}},j)$. Each path $P$ that traverses $e_+$ exactly once can be decomposed as $P=(P_1,e_-,e_+,P_2)$ so that $\hP_1=(P_1,\hat e_-)$ is a path from $b_{i_{\bar p}}$ to $b_{j_u}$ in $\hN$ and $\hP_2=(\hat e_+, P_2)$ is a path from $b_{i_u}$ to $b_j$ in $\hN$. Define $\hP=(\hP_1,e,\hP_2)$, where $e$ is the edge between $j_u$ and $i_u$ belonging to the convex hull of the boundary vertices of $\hN$ (see Figure~\ref{fig:recalc}). Clearly, $c(C_P)=c(C_{\hP})$. 

\begin{figure}[ht]
\begin{center}
\includegraphics[height=12cm]{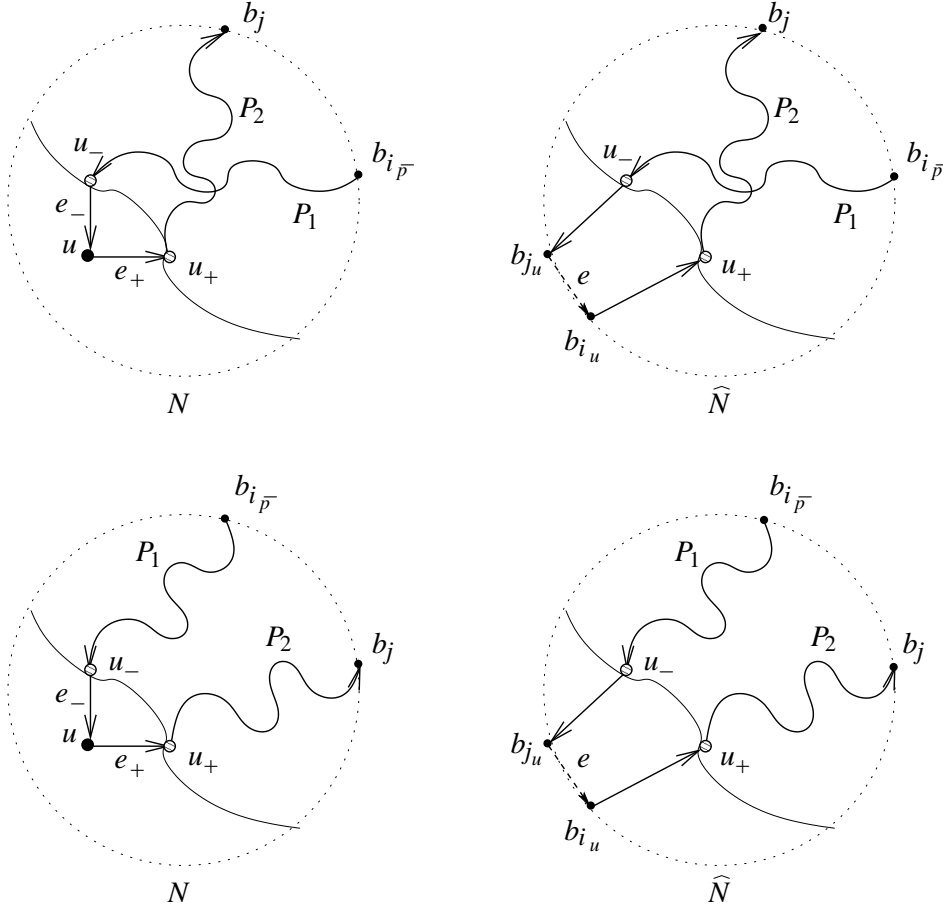}
\caption{To the proof of Lemma~\ref{recalc}} \label{fig:recalc}
\end{center}
\end{figure}

Assume first that 
\[i_{\bar p}\prec j\preceq i_p-1\prec j_u\prec i_u\prec i_p+1, 
\]
which corresponds to the upper part of Figure~\ref{fig:recalc}. Then 
\[
c(C_{\hP_1})=c_{\hP_1}+c_{i_u}+c', 
\]
where $c_{\hP_1}$ is the contribution of all vertices of $\hP_1$, including $b_{i_{\bar p}}$ and $b_{j_u}$, $c_{i_u}$ is the contribution of two consecutive edges of the convex hull of boundary vertices of $\hN$ calculated at $b_{i_u}$, and $c'$ is the total contribution calculated at the vertices of the convex hull lying between $b_{i_u}$ and $b_{i_{\bar p}}$. Similarly, 
\[
c(C_{\hP_2})=c_{\hP_2}+c_{j_u}+c'', 
\]
where $c_{\hP_2}$ is the contribution of all vertices of $\hP_2$, including $b_{i_u}$ and $b_{j}$, $c_{j_u}$ is the contribution of two consecutive edges of the convex hull of boundary vertices of $\hN$ calculated at $b_{j_u}$, and $c''$ is the total contribution calculated at the vertices of the convex hull lying between $b_{j}$ and $b_{j_{u}}$. Finally, 
\[
c(C_{\hP})=c_{\hP_1}+c_{\hP_2}+c''+c_{j_u}+
c_{i_u}+c', 
\]
and so $c(C_{\hP})=c(C_{\hP_1})+c(C_{\hP_2})$. Therefore, in this case $w_P=-w_{P_1}w_{e_-}w_{e_+}w_{P_2}$, and the contribution of all such paths to $M({i_{\bar p}},j)$ equals $-x''\hM({i_{\bar p}},{j_u})\hM({i_u},j)$. Each additional traversing of the edge $e_+$ results in multiplying this expression by $-x''\hM({i_u},{j_u})$; the proof of this fact is similar to the proof of Proposition~\ref{sfree}. Summing up we get the second formula with the $-$ sign, as desired.

Assume now that 
\[
i_p-1\prec j_u\prec i_u\prec i_p+1\preceq j \prec i_{\bar p}, 
\]
which corresponds to the lower part of Figure~\ref{fig:recalc}. Then 
\[
c(C_{\hP_1})=c_{\hP_1}+c_{i_u}+c'+c_j+c''', 
\]
where $c_{\hP_1}$ and $c_{i_u}$ are as in the previous case, $c_j$ is the contribution of two consecutive edges of the convex hull of boundary vertices of $\hN$ calculated at $b_j$, $c'$ is the total contribution calculated at the vertices of the convex hull lying between $b_{i_u}$ and $b_{j}$, and $c'''$ is the total contribution calculated at the vertices of the convex hull lying between $b_j$ and $b_{i_{\bar p}}$. Similarly, 
\[
c(C_{\hP_2})=c_{\hP_2}+c'''+c_{i_{\bar p}}+c''+c_{j_u}, 
\]
where $c_{\hP_2}$ and $c_{j_u}$ are as in the previous case, $c_{i_{\bar p}}$ is the contribution of two consecutive edges of the convex hull of boundary vertices of $\hN$ calculated at $b_{i_{\bar p}}$, and $c''$ is the total contribution calculated at the vertices of the convex hull lying between $b_{i_{\bar p}}$ and $b_{j_u}$. Finally, $c(C_{\hP})=c_{\hP_1}+c_{\hP_2}+c'''$, and so $c(C_{\hP})=c(C_{\hP_1})+c(C_{\hP_2})+1$, since 
\[
c'+c''+c'''+c_{j_u}+c_{i_u}+ c_j+c_{i_{\bar p}}=0. 
\]
Therefore, in this case $w_P=w_{P_1}w_{e_-}w_{e_+}w_{P_2}$, and the contribution of all such paths to $M({i_{\bar p}},j)$ equals $x''\hM({i_{\bar p}},{j_u})\hM({i_u},j)$. Each additional traversing of $e_+$ once again results in multiplying this expression by $-x''\hM({i_u},{j_u})$. Summing up we get the second formula with the $+$ sign, as desired.

To treat the remaining two cases one makes use of Remark~\ref{dualdir} 
and applies the same reasoning. 
\end{proof}

Since $\hN$ has $r-1$ inner vertices, 
$\hM(i,j)$ satisfy relations similar to~(\ref{twops}) with $N$ replaced by $\hN$, $I$ replaced by $I\setminus i_p\cup i_v$ and $J$ replaced by $J\cup j_v$.
Relations~(\ref{twops}) follow from this fact and Lemma~\ref{recalc} via simple though tedious computations.

For example, let $i_{\bar p}\prec j_{\bar q}\preceq i_p-1\prec j_u\prec i_u\prec i_p+1$. Then the left hand side of~(\ref{twops}) is given by
\begin{align*}
\{M_{pq},M_{p\bar q}\}_N
&=\left\{\frac{x'\hM({i_u},j_q)}{1+x''\hM({i_u},{j_u})},
\hM({i_{\bar p}},j_{\bar q})
-\frac{x''\hM({i_{\bar p}},{j_u})\hM({i_u},j_{\bar q})}
{1+x''\hM({i_u},{j_u})}\right\}_N\\
&=\frac{x'}{1+x''\hM({i_u},{j_u})}\{\hM(i_u,j_q),
\hM(i_{\bar p},j_{\bar q})\}_\hN\\
&-\frac{x'x''\hM(i_u,j_q)}{(1+x''\hM({i_u},{j_u}))^2}
\{\hM(i_u,j_u),\hM(i_{\bar p},j_{\bar q})\}_\hN\\
&-\frac{\hM({i_u},j_q)\hM({i_{\bar p}},j_{u})\hM(i_{u},j_{\bar q})}
{(1+x''\hM({i_u},{j_u}))^2}\{x',x''\}_u\\
&-\frac{x'x''\hM(i_u,j_{\bar q})}{(1+x''\hM({i_u},{j_u}))^2}
\{\hM(i_u,j_q),\hM(i_{\bar p},j_u)\}_\hN\\
&-\frac{x'x''\hM(i_{\bar p},j_{u})}{(1+x''\hM({i_u},{j_u}))^2}
\{\hM(i_u,j_q),\hM(i_{u},j_{\bar q})\}_\hN\\
&+\frac{x'(x'')^2\hM(i_{u},j_{\bar q})\hM(i_{u},j_{q})}
{(1+x''\hM({i_u},{j_u}))^3}\{\hM(i_u,j_u),\hM(i_{\bar p},j_{u})\}_\hN\\
&+\frac{x'(x'')^2\hM(i_{u},j_{q})\hM(i_{\bar p},j_{u})}
{(1+x''\hM({i_u},{j_u}))^3}\{\hM(i_u,j_u),\hM(i_{u},j_{\bar q})\}_\hN\\
&+\frac{x'(x'')^2\hM(i_{\bar p},j_{u})\hM(i_{u},j_{\bar q})}
{(1+x''\hM({i_u},{j_u}))^3}\{\hM(i_u,j_q),\hM(i_{u},j_{u})\}_\hN\\
&+\frac{x''\hM(i_{\bar p},j_{u})\hM(i_{u},j_{\bar q})
\hM(i_u,j_q),\hM(i_{u},j_{u})}
{(1+x''\hM({i_u},{j_u}))^3}\{x',x''\}_u,
\end{align*}
and the right hand side is given by
\begin{align*}
(&\alpha-\beta)s_=(i_p,j_q,i_{\bar p},j_{\bar q})
\frac{x'\hM({i_u},j_{\bar q})}{1+x''\hM({i_u},{j_u})}\left(
\hM({i_{\bar p}},j_{q})
-\frac{x''\hM({i_{\bar p}},{j_u})\hM({i_u},j_{q})}
{1+x''\hM({i_u},{j_u})}\right)\\
&+(\alpha+\beta)s_\times(i_p,j_q,i_{\bar p},j_{\bar q})
\frac{x'\hM({i_u},j_q)}{1+x''\hM({i_u},{j_u})}\left(
\hM({i_{\bar p}},j_{\bar q})
-\frac{x''\hM({i_{\bar p}},{j_u})\hM({i_u},j_{\bar q})}
{1+x''\hM({i_u},{j_u})}\right).
\end{align*}
Treating $\hM(\cdot,\cdot)$, $\alpha-\beta$ and $\alpha +\beta$ as independent variables and equating coefficients of the same monomials in the above two expressions we arrive to the following identities:
\begin{align*}
&s_=(i_u,j_q,i_{\bar p},j_{\bar q})=s_=(i_p,j_q,i_{\bar p},j_{\bar q}),\qquad
s_\times(i_u,j_q,i_{\bar p},j_{\bar q})=
s_\times(i_p,j_q,i_{\bar p},j_{\bar q}),\\
&s_=(i_u,j_u,i_{\bar p},j_{\bar q})+s_=(i_u,j_q,i_{\bar p},j_{u})
+s_=(i_u,j_q,i_{u},j_{\bar q})-1/2=s_=(i_p,j_q,i_{\bar p},j_{\bar q}),\\
&s_\times(i_u,j_u,i_{\bar p},j_{\bar q})+s_\times(i_u,j_q,i_{\bar p},j_{u})
+s_\times(i_u,j_q,i_{u},j_{\bar q})+1/2=
s_\times(i_p,j_q,i_{\bar p},j_{\bar q}),\\
&s_=(i_u,j_u,i_{\bar p},j_{u})+s_=(i_u,j_u,i_{u},j_{\bar q})+
s_=(i_u,j_q,i_{u},j_{u})-1/2=0,\\
&s_\times(i_u,j_u,i_{\bar p},j_{u})+s_\times(i_u,j_u,i_{u},j_{\bar q})+
s_\times(i_u,j_q,i_{u},j_{u})+1/2=0.
\end{align*}
The latter can be checked easily by considering separately the following cases: $i_p\prec j_q\prec i_{\bar p}$; $i_{\bar p}\prec j_q\prec j_{\bar q}$;
$j_q=j_{\bar q}$; $j_{\bar q}\prec j_q\prec i_p$. In each one of these cases all the functions involved in the above identities take constant values.
\end{proof}

\section{Grassmannian boundary measurement map and induced Poisson structures on $G_k(n)$}\label{PSongr}

\subsection{Pl\"ucker coordinates and Poisson brackets}
\label{recallGSV1}

Let us recall the  construction of Poisson brackets on the Grassmannian from \cite{GSV1}.
Let $G_k(n)$ be the Grassmannian of $k$-dimensional subspaces in $\R^n$.
Given  a $k$-element  subset  $I$  of $[1,n]=\{1,\dots,n\}$,
the  Pl\"ucker coordinate\/ $x_I$ 
is a function on the set of $k\times n$ matrices which is equal to the value of the minor formed by the columns of the matrix indexed by the elements of $I$.
 In what follows, we use
notation
\begin{equation}
I(i_\alpha\to l) = \{i_1, \ldots, i_{\alpha-1}, l ,i_{\alpha+1},
\ldots, i_k \}
\nonumber
\end{equation}
for $\alpha\in [1, k]$ and $l\in [1,n]$.

Let  $G^0_k(n)$
be the open cell in $G_k(n)$ characterized by
non-vanishing of the Pl\"ucker coordinate $x_{[1, k]}$. 
Elements of $G^0_k(n)$ are parametrized by $k\times(n-k)$ matrices
in the following way: if $W\in SL_n$ admits a factorization into
block-triangular matrices
\begin{equation}
W= \begin{pmatrix} W_{11} & 0\\ W_{12} & W_{22} \end{pmatrix}
\begin{pmatrix} \one_k &  Y \\ 0 & \one_{n-k} \end{pmatrix},
\label{defY}
\end{equation}
then $Y=Y(W)$ represents an element of the cell $G^0_k(n)$.

It is easy to check that Pl\"ucker coordinates $x_I$,
$I=\{i_1,\ldots,i_k\ : 1\ls i_1 < \cdots < i_k \ls n\}$ of an element of $G^0_k(n)$ represented by $[\one_k\ Y]$ and
minors $Y_{\alpha_1, \ldots, \alpha_l}^{\beta_1, \ldots, \beta_l}=
\det(y_{\alpha_i,\beta_j})_{i,j=1}^l$ of $Y$ are related  via
$$
Y_{\alpha_1, \ldots, \alpha_l}^{\beta_1, \ldots, \beta_l} = (-1)^{k l
-l(l-1)/2 - (\alpha_1 + \cdots + \alpha_l)}
\frac{x_{([1,k]\setminus\{\alpha_1, \ldots, \alpha_l\}) \cup
\{\beta_1+k, \ldots, \beta_l+k\}}}{x_{[1,k]}}.
$$
Note that, if the row index set $\{\alpha_1, \ldots, \alpha_l\}$ in
the above formula is contiguous then the sign in the right hand
side can be expressed as $(-1)^{(k-\alpha_l) l}$.

A Poisson bracket $\Poi$ on
$G^0_k(n)$ is defined via the relation
$$
\{f_1\circ Y,f_2\circ Y\}_{\Mat_n}=\{ f_1 , f_2\}\circ Y\ .
$$

In terms of matrix elements $y_{ij}$ of $Y$, this bracket  looks
as follows:
\begin{equation}\label{brackY}
 2 \{y_{ij}, y_{\alpha\beta}\}= (\s (\alpha - i)- \s (\beta -j))
y_{i\beta} y_{\alpha j}.
\end{equation}

\subsection{Recovering the Sklyanin bracket on $\Mat_k$}
Let us take a closer look at the 2-parameter family of Poisson brackets obtained in Theorem~\ref{PSME} in the case when vertices $b_1,\dots,b_k$ on the boundary of the disk are  sources and vertices $b_{k+1},\dots,b_n$ are sinks, that is, when $I=[1,k]$ and $J=[k+1,n]$. To simplify notation, in this situation we will write $\{\cdot,\cdot\}_{k,m}$ and $\Net_{k,m}$ instead of $\{\cdot,\cdot\}_{[1,k],[k+1,n]}$ and $\Net_{[1,k],[k+1,n]}$. Therefore,
formula~(\ref{psme}) can be re-written as
\begin{multline}\label{Mbrack12}
2\{M_{ij}, M_{\bar\imath\bar\jmath}\}_{k,m}= (\alpha-\beta)\left (\s (\bar\imath - i)- \s (\bar\jmath -j)\right )
M_{i\bar\jmath} M_{\bar\imath j}\\
+ (\alpha+\beta) \left (\s (\bar\imath - i)+ \s (\bar\jmath -j)\right )
M_{ij} M_{\bar\imath\bar\jmath},
\end{multline}
where $M_{ij}$ corresponds to the boundary measurement between $b_i$ and $b_{j+k}$.
The first term in the equation above coincides  (up to a multiple) with (\ref{brackY}).
This suggests that it makes sense to investigate  Poisson properties of the boundary measurement map viewed as a map into $G_k(n)$, which will be the goal of this section.

First, however, we will go back to the example, considered in Sect.~\ref{concat}, where we associated with a network
$N\in \Net_{k,k}$  a matrix $A_N = M_N W_0\in \Mat_k$. Written in terms of matrix entries $A_{ij}$ of $A$, bracket~(\ref{Mbrack12}) becomes
\begin{multline}\label{Abrack12}
2\{A_{ij}, A_{\bar\imath\bar\jmath}\}_{k,m}= (\alpha-\beta)\left (\s (\bar\imath - i)+ \s (\bar\jmath -j)\right )
A_{i\bar\jmath} A_{\bar\imath j}\\ 
+ (\alpha+\beta) \left (\s (\bar\imath - i)- \s (\bar\jmath -j)\right )
A_{ij} A_{\bar\imath\bar\jmath}.
\end{multline}
If $A_{N_1}, A_{N_2}$ are matrices that correspond to networks $N_1, N_2 \in \Net_{k,k}$ then their product
$A_{N_1} A_{N_2}$ corresponds to the concatenation of $N_1$ and $N_2$. This fact combined with Theorem~\ref{PSM} 
implies that the bracket (\ref{Abrack12}) possesses the Poisson-Lie property. 

In fact, (\ref{Abrack12}) is the Sklyanin
bracket  ~(\ref{sklyaSLn}) on $\Mat_k$ 
associated with a deformation of the standard R-matrix~(\ref{standardR}). Indeed,
it is known (see \cite{r-sts}) that if $R_0$ is the standard R-matrix,
$S$ is any linear operator on the set of diagonal matrices 
that is skew-symmetric w.r.t. the trace-form, and $\pi_0$ is the natural projection onto a subspace of diagonal matrices, then for any scalar $c_1$, $c_2$ the linear combination $ c_1 R_0 + c_2 S \pi_0$ satisfies MCYBE  (\ref{MCYBE}) and thus gives rise to a Sklyanin Poisson-Lie bracket.

Define $S$ by
\begin{equation*}
S(e_{jj}) = \sum_{i=1}^k \s(j-i) e_{ii},\quad j=1,\ldots, k,
\end{equation*}
and put
\begin{equation*}
R_{\alpha,\beta}=\frac{\alpha-\beta}{2} R_0 + \frac{\alpha+\beta}{2} S \pi_0.
\end{equation*}
 Substituting coordinate functions $A_{ij}$, $A_{i'j'}$ into expression~(\ref{sklyaSLn}) for the  Sklyanin Poisson-Lie bracket 
associated  with the R-matrix $R_{\alpha,\beta}$, we recover equation (\ref{Abrack12}).

To summarize, we obtained
\begin{theorem} \label{PSGL}
For any network $N \in \Net_{k,k}$ and any choice of parameters $\alpha_{ij}$, $\beta_{ij}$ in~\eqref{6parw},~\eqref{6parb},
the map $A_N: \RS^{|E|} \to \Mat_k$ is Poisson with respect to the Sklyanin bracket
~\eqref{sklyaSLn} 
associated with the R-matrix $R_{\alpha,\beta}$, where $\alpha$ and $\beta$ satisfy relations~\eqref{condw} and~\eqref{condb}.
\end{theorem}  

\subsection{Induced Poisson structures on $G_k(n)$}\label{indpoisson} 

Let $N$ be a network with the sources $b_i$, $i\in I$ and 
sinks $b_j$, $j\in J$.
Following \cite{Postnikov}, we are going to interpret the boundary measurement
map as a map into the Grassmannian $G_k(n)$. To this end, we extend $M_N$ to a $k\times n$ matrix $\wX_N$ as follows: 

(i) the $k\times k$ submatrix of $\wX_N$ formed
by $k$ columns indexed by $I$ is the identity matrix $\one_k$; 

(ii) for $p\in[1,k]$ and $j=j_q\in J$, the 
$(p,j)$-entry of $\wX_N$ is 
\[
m^I_{pj}= (-1)^{s(p,j)} M_{p q}, 
\]
where $s(p,j)$ is the number of elements in $I$ lying strictly between $\min\{i_p,j\}$ and $\max\{i_p,j\}$ in the linear ordering; note that the sign is selected in such a way that the minor
$(\wX_N)_{[1,k]}^{I(i_p\to j)}$ coincides with $M_{p q}$. 

\begin{figure}[ht]
\begin{center}
\includegraphics[height=4cm]{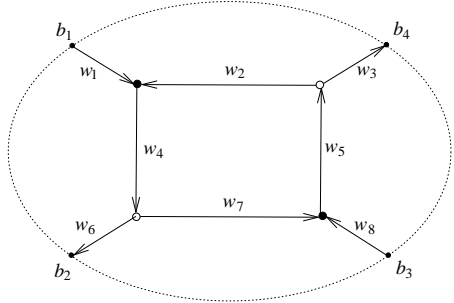}
\caption{To the definition of the Grassmannian boundary measurement map} \label{fig:GrBMM}
\end{center}
\end{figure}

We will view $\wX_N$ as a matrix representative of an element $X_N \in G_k(n)$. 
The corresponding rational map $X_N\: \CC_N\to G_k(n)$ is called the {\it Grassmannian boundary measurement map}. For example, the network presented in Figure~\ref{fig:GrBMM} defines a map of $\RS^8$ to $G_2(4)$ given by the matrix
\[
\begin{pmatrix}
&1 &\dfrac{w_1w_4w_6}{1+w_2w_4w_5w_7}           &0  &-\dfrac{w_1w_3w_4w_5w_7}{1+w_2w_4w_5w_7}\\
\\
&0 &\dfrac{w_2w_4w_5w_6w_8}{1+w_2w_4w_5w_7}     &1
&\dfrac{w_3w_5w_8}{1+w_2w_4w_5w_7}
\end{pmatrix}.
\]

Clearly, $X_N$ belongs to the cell 
$G_k^I(n)=\{ X \in  G_k(n)\ : \ x_I \ne 0\}$. Therefore, we can regard matrix entries $m^I_{pj}$ 
as coordinate functions on $G_k^I(n)$ and 
rewrite~(\ref{psme}) as follows: 
\begin{equation}\label{mbrack12}
\{m^I_{pj}, m^I_{\bar p\bar \jmath}\}_{I,J}=(\alpha-\beta)s_=(i_p, j, i_{\bar p},\bar\jmath) m^I_{p\bar\jmath} m^I_{\bar pj} 
+  (\alpha+\beta)s_\times(i_p, j, i_{\bar p},\bar\jmath) m^I_{pj}m^I_{\bar p\bar\jmath}.
\end{equation}

\begin{remark}{\rm
If $I=[1,k]$, then the 2-parametric family~(\ref{mbrack12}) of Poisson brackets on $G_k^I(n)=G_k^0(n)$ is defined by~(\ref{Mbrack12}). A computation identical to the one presented in \cite{GSV1} (see Sect. \ref{recallGSV1}) shows that~(\ref{Mbrack12}) is the pushforward to $G_k^0(n)$
of the Sklyanin Poisson-Lie bracket~(\ref{sklyaSLn}) on $\Mat_n$ associated with the R-matrix $R_{\alpha,\beta}$. }
\end{remark}

The following result says that the families $\{\cdot,\cdot\}_{I,J}$ on different cells $G^I_k(n)$ can be glued together to form the unique 2-parametric family of Poisson brackets on $G_k(n)$ that makes all maps $X_N$ Poisson.

\begin{theorem}\label{PSGr}
{\rm(i)} For any choice of parameters $\alpha$ and $\beta$ there exists a unique Poisson bracket $\P_{\alpha,\beta}$ on $G_k(n)$ such that for any network $N$ with $k$ sources, $n-k$ sinks and weights defined by~\eqref{connect},  the map $X_N\: \RS^{|E|}\to G_k(n)$ is Poisson provided the parameters $\alpha_{ij}$ and $\beta_{ij}$ defining the bracket $\{\cdot,\cdot\}_N$ on $\RS^{|E|}$ satisfy relations~\eqref{condw} and~\eqref{condb}.

{\rm(ii)} For any $I\subset [1,n]$, $|I|=k$, and $J=[1,n]\setminus I$, the restriction of $\P_{\alpha,\beta}$ to the cell $G^I_k(n)$ coincides with the bracket $\{\cdot,\cdot\}_{I,J}$ given by~\eqref{mbrack12} in coordinates $m^I_{pj}$.
\end{theorem}

\begin{proof} Both statements follow immediately from Theorems~\ref{PSM} and~\ref{PSME}, the fact that $G_k^0(n)$ is an open dense subset in $G_k(n)$, and the following proposition.

\begin{proposition}\label{Coincbr}
For any $i_p,i_{\bar p}\in I$ and any $j,\bar\jmath\in J$,
\begin{equation}\label{coincbr}
\{ m^I_{p j}, m^I_{{\bar p} \bar\jmath}\}_{k,n-k}=\{ m^I_{p j}, m^I_{{\bar p} \bar\jmath}\}_{I,J}. 
\end{equation}
\end{proposition}

\begin{proof}
It is convenient to rewrite~(\ref{mbrack12}) as 
\[
\{m^I_{pj}, m^I_{\bar p\bar \jmath}\}_{I,J}=\{m^I_{pj}, m^I_{\bar p\bar \jmath}\}^1_{I,J}+\{m^I_{pj}, m^I_{\bar p\bar \jmath}\}^2_{I,J} 
\]
and to treat the brackets $\{\cdot,\cdot\}^1_{I,J}$ and  $\{\cdot,\cdot\}^2_{I,J}$ separately for a suitable choice of parameters $\alpha$ and $\beta$. 

Let $I=[1,k]$ and denote $m_{ij}^{[1,k]}$ simply by $m_{ij}$. Then  two independent brackets are given by
\begin{equation}\label{firstbrack}
  \{m_{ij}, m_{\bar\imath\bar\jmath}\}^1_{k,n-k}= \left (\s (\bar\imath - i)- \s (\bar\jmath -j)\right )
m_{i\bar\jmath} m_{\bar\imath j}.
\end{equation}
(for $\alpha-\beta=2$)
that coincides (up to a multiple) with the Poisson 
bracket  
(\ref{brackY}) on $G^0_k(n)$ and
\begin{equation}\label{secondbrack}
  \{m_{ij}, m_{\bar\imath\bar\jmath}\}^2_{k,n-k}= \left (\s (\bar\imath - i)+ \s (\bar\jmath -j)\right )
m_{ij} m_{\bar\imath \bar\jmath}
\end{equation}
(for $\alpha+\beta=2$).

Let us start with the first of the two brackets above. For a generic element $X$ in  $G_k(n)$ represented by a matrix $\wX=\wX_N$,
consider any two minors, $\wX_\mu^\nu$ and $\wX_{\mu'}^{\nu'}$,  
with row sets $\mu=(\mu_1,\ldots,\mu_l)$,  $\mu'=(\mu'_1,\ldots,\mu'_{l'})$, $\mu,\mu' \subseteq [1,k]$ and column sets $\nu=(\nu_1,\ldots,\nu_l)$, $\nu'=(\nu'_1,\ldots,\nu'_{l'})$, $\nu,\nu'\subseteq [k+1,n]$.
Using considerations similar to those in the proof of Lemma 3.2  in \cite{GSV1},
we obtain from~(\ref{firstbrack})
\begin{multline}\label{firstbrackminors}
\{\wX_\mu^\nu, \wX_{\mu'}^{\nu'} \}^1_{k,n-k} =
\sum_{p=1}^l  \sum_{q=1}^{l'} \left ( 
\s (\mu'_p - \mu_q) \wX_{\mu(\mu_p\to \mu'_q)}^\nu \wX_{\mu'(\mu'_q\to \mu_p)}^{\nu'} \right.\\
-  \left.\s (\nu'_p - \nu_q) \wX_\mu^{\nu(\nu_p\to \nu'_q)}  \wX_{\mu'}^{\nu'(\nu'_q\to \nu_p)} \right ).
\end{multline}

 Pl\"ucker coordinates $x_I$ (for $I\subset [1,n]$, $|I|=k$)  and minors $\wX_{\mu}^{\nu}$ are related  via
$$
\wX_{\mu}^{\nu} = (-1)^{
l(l-1)/2}
\frac{x_{([1,k]\setminus \mu) \cup
\{\nu\}}}{x_{[1,k]}}.
$$
Denote
$$
a_I = \frac{x_{I}}{x_{[1,k]}} = (-1)^{
l(l-1)/2} \wX_{[1,k]\setminus I}^{I\setminus[1,k]},
$$
where $l=\left  |[1,k]\setminus I \right |$.
Then~(\ref{firstbrackminors}) gives rise to Poisson relations
\begin{eqnarray*}
\{a_I, a_{I'} \}^1_{k,n-k} &=
\sum_{i\in I\cap [1,k] } \sum_{i'\in I'\cap [1,k]}
\s (i - i')a_{I(i\to i')} a_{I'(i'\to i)}\\
&+\sum_{i\in I\setminus [1,k] } \sum_{i'\in I'\setminus [1,k] }
\s (i - i')  a_{I(i\to i')} a_{I'(i'\to i)}\\ \label{firstbrackPlu}
&=\sum_{i\in I} \sum_{i'\in I'} \varepsilon_{i i'}  
a_{I(i\to i')} a_{I'(i'\to i)},
\end{eqnarray*}
where
$$
 \varepsilon_{i i'} = \left \{ \begin{array}{cc} 
0& \mbox{if} \quad i\le k < \bar\imath \quad \mbox{or} \quad \bar\imath\le k < i,\\
 \s (i-i') & \mbox{otherwise}. \end{array} \right .
 $$

Let us fix a $k$-element index set $I=\{1\le i_1< \ldots < i_k \le n\}$ and 
use~(\ref{firstbrackPlu}) to  compute, for any $p,\bar p\in[1,k]$ and $j, \bar\jmath \notin I$, the bracket
$\{m_{pj}^I,m_{\bar p\bar\jmath}^I\}^1_{k,n-k}$. Taking into account that 
\[
m^I_{pj}= \frac{x_{I(i_p\to j)}}{x_I}\quad\text{for any $p\in [1,k]$, $j \in J$},
\]
we get
\begin{align*}
\{ m^I_{p j}, m^I_{{\bar p} \bar\jmath}\}^1_{k,n-k} &= 
\dfrac{1}{a_{I}^2} 
\bigg(\{a_{I(i_p\to j)}, a_{I(i_{\bar p}\to \bar\jmath)} \}^1_{k,n-k} - \dfrac{a_{I(i_p\to j)}}{a_{I}}  
\{a_{I}, a_{I(i_{\bar p}\to \bar\jmath)} \}^1_{k,n-k}\\  
  &\qquad\qquad \dfrac{a_{I(i_{\bar p}\to \bar\jmath)}}{a_I}  \{a_{I(i_p\to j)}, a_{I} \}^1_{k,n-k}  \bigg )\\
&= \dfrac{1}{a_{I}^2} \big( 
 (\varepsilon_{j\bar\jmath} + \varepsilon_{i_{\bar p}  i_p})  a_{I(i_p\to \bar\jmath)} a_{I(i_{\bar p}\to j)} 
+(\varepsilon_{j i_p} + \varepsilon_{i_{\bar p}\bar\jmath})a_I a_{I(i_p\to j, i_{\bar p} \to \bar\jmath)}  \\
 &\qquad\qquad - (\varepsilon_{j i_p} + \varepsilon_{i_{\bar p}\bar\jmath})  
a_{I(i_p\to j)} a_{I(i_{\bar p}\to \bar\jmath)} \big) \\
&=\dfrac{1}{a_I^2}  ( \varepsilon_{j\bar\jmath} + \varepsilon_{i_{\bar p}  i_p} - \varepsilon_{j i_p} -  \varepsilon_{i_{\bar p}  \bar\jmath}) a_{I(i_p\to \bar\jmath)} a_{I(i_{\bar p}\to j)}\\
& =(\varepsilon_{j\bar\jmath}+\varepsilon_{i_{\bar p}i_p}-\varepsilon_{j i_p} -\varepsilon_{i_{\bar p}\bar\jmath})  m^I_{p \bar\jmath} m^I_{{\bar p} j},
\end{align*}
 where in the third step we have used the short Pl\"ucker relation.

Relation~(\ref{coincbr}) for the first bracket follows from

\begin{lemma}
For any $i_p,i_{\bar p}\in I$ and any $j,\bar\jmath\notin I$,
\begin{equation}\label{newspar}
  \varepsilon_{j\bar\jmath} + \varepsilon_{i_{\bar p}  i_p}-\varepsilon_{j i_p}-\varepsilon_{i_{\bar p}\bar\jmath}=
2s_=(i_p,j,i_{\bar p},\bar\jmath).
\end{equation}
\end{lemma}

\begin{proof}
Denote the left hand side of~(\ref{newspar}) by $\varepsilon^1(i_p,j,\bar\jmath,i_{\bar p})$. Let us prove  
that this expression depends only on the counterclockwise order of the numbers $i_p$, $j$, $\bar\jmath$, $i_{\bar p}$ and does not depend on the numbers themselves. 

Assume first that all four numbers are distinct. In this case $\varepsilon^1(i_p,j,\bar\jmath,i_{\bar p})$ is invariant with respect to the cyclic shift of the variables, hence it suffices to verify the identity 
$\varepsilon^1(i_p,j,\bar\jmath,i_{\bar p}) =\varepsilon^1(i_p,j-1\mod n,\bar\jmath,i_{\bar p})$ provided the counterclockwise orders of $i_p$, $j$, $\bar\jmath$, $i_{\bar p}$ and $i_p$, $j-1\mod n$, $\bar\jmath$, $i_{\bar p}$ coincide. This is trivial unless $j=1$ or $j=k$, since all four summands retain their values. In the remaining cases the second and the fourth summands retain their values, and we have to check the identities 
$\varepsilon_{1\bar\jmath}-\varepsilon_{1 i_p}=\varepsilon_{n\bar\jmath}-\varepsilon_{n i_p}$ and
$\varepsilon_{k\bar\jmath}-\varepsilon_{k i_p}=\varepsilon_{k-1,\bar\jmath}-\varepsilon_{k-1, i_p}$. The first of them follows from the identity $\varepsilon_{ni}=\varepsilon_{1i}+1$ for $i\ne 1,n$, and the second one, from the identity $\varepsilon_{k-1,i}=\varepsilon_{ki}+1$ for $i\ne k-1,k$.

If $j=\bar\jmath$ or $i_p=i_{\bar p}$ (other coincidences are impossible since $i_p,i_{\bar p}\in I$, $j,\bar\jmath\notin I$), $\varepsilon^1$ degenerates to a function of three variables, which is again invariant with respect to cyclic shifts, therefore all the above argument applies as well. 

To obtain~(\ref{newspar}) it remains to check that  $\varepsilon^1(i_p,1,k,i_{\bar p})=2s_=(i_p,1,i_{\bar p},k)$, which can be done separately in all the cases mentioned in~(\ref{spar}).
\end{proof}

Next, let us turn to the Poisson bracket~(\ref{secondbrack}). For any two monomials in matrix entries of $\wX_N$, 
 $\mathbf{x}=\prod_{p=1}^l \wX_{\mu_p \nu_p}$, and  $\mathbf{x}'=\prod_{q=1}^{l'} \wX_{\mu'_q \nu'_q}$, we have
 $$
 \{ \mathbf{x}, \mathbf{x}'\}^2_{k,n-k} = \mathbf{x}\  \mathbf{x}' \sum_{p=1}^l \sum_{q=1}^{l'} 
 \left ( \s (\mu'_p - \mu_q) + \s (\nu'_p - \nu_q)\right).
 $$
 Observe that the double sum above is invariant under any permutation of indices within sets
 $\mu=(\mu_1,\ldots,\mu_l)$, $\mu'=(\mu'_1,\ldots,\mu'_{l'})$, $\nu=(\nu_1,\ldots,\nu_l)$, 
$\nu'=(\nu'_1,\ldots,\nu'_{l'})$. This means that for
 minors  $\wX_\mu^\nu, \wX_{\mu'}^{\nu'}$  we have
\begin{equation*}
\{\wX_\mu^\nu, \wX_{\mu'}^{\nu'} \}^2_{k,n-k} = \wX_\mu^\nu \wX_{\mu'}^{\nu'}\sum_{p=1}^l \sum_{q=1}^{l'} 
 \left ( \s (\mu'_p - \mu_q) + \s (\nu'_p - \nu_q)\right ), 
\end{equation*}
and the resulting Poisson relations for functions $a_I$ are
\begin{align*}
\{a_I, a_{I'} \}^2_{k,n-k} &= a_I a_{I'} \left(\sum_{i\in [1,k] \setminus I} \sum_{i'\in  [1,k] \setminus I'} \s (i - i')\right.\\
&\qquad\qquad\qquad +\left.
\sum_{i\in I\setminus [1,k] } \sum_{i'\in I'\setminus [1,k] }\s (i - i') 
\right)\\ 
&= a_I a_{I'}\sum_{i\in I} \sum_{i'\in I'} \varepsilon_{i i'},
\end{align*}
where $\varepsilon_{ii'}$ has the same meaning as above.
 These relations imply
 $$
  \{ m^I_{i_p j}, m^I_{i_{\bar p} \bar\jmath}\}^2_{k,n-k} = m^I_{i_p j} m^I_{i_{\bar p} \bar\jmath} 
\left(\varepsilon_{i_p i_{\bar p}} - \varepsilon_{i_p \bar\jmath} - \varepsilon_{q i_{\bar p}} + 
\varepsilon_{j\bar\jmath}   \right ).
 $$

Relation~(\ref{coincbr}) for the second bracket follows now from

\begin{lemma}
For any $i_p,i_{\bar p}\in I$ and any $j,\bar\jmath\notin I$,
\begin{equation}\label{newscross}
 \varepsilon_{i_p i_{\bar p}} - \varepsilon_{i_p \bar\jmath} - \varepsilon_{q i_{\bar p}} + 
\varepsilon_{j\bar\jmath}  = 2s_\times(i_p,j,i_{\bar p},\bar\jmath).
\end{equation}
\end{lemma}

\begin{proof} The proof of~(\ref{newscross}) is similar to the proof of~(\ref{newspar}).
\end{proof}

This proves Proposition~\ref{Coincbr}, and hence Theorem~\ref{PSGr}.
\end{proof}

\end{proof}

\subsection{$GL_n$-action on the Grassmannian via networks}
In this subsection we interpret the natural action of $GL_n$ on $G_k(n)$ in terms of planar networks. 

First, note
that {\em any} element of $GL_n$ can be represented by a planar network built by concatenation
from building blocks (elementary networks)  described in Fig.~\ref{fig:factor}. To see this, one needs to observe that an elementary
transposition matrix $S_i = \one - e_{ii} - e_{i+1,i+1} + e_{i,i+1} + e_{i+1,i}$ can be factored as
\[
S_i= (\one - 2 e_{i+1,i+1})  E_{i+1}^-(-1) E_{i+1}^+(1) E_{i+1}^-(-1), 
\]
which implies that any permutation matrix can be represented via concatenation of elementary networks. Consequently, any elementary 
triangular matrices $\one + u e_{i j}$,  $\one + l e_{ji}$,  $i< j$, can be factored as
\[
 \one + u e_{i j}= W_{ij} E_{i+1}^+(u) W^{-1}_{ij},\qquad 
 \one + l e_{ji }= W^{-1}_{ij} E_{i+1}^-(l) W_{ij},
\] 
 where $W_{ij}$ is the permutation
 matrix that corresponds to the permutation
 \[ 
 (1)\cdots (i)  (i+1 \ldots j) (j+1)\cdots (n). 
 \]
 The claim then follows from the Bruhat decomposition and constructions presented in Section~\ref{concat}.

 Consider now a network $N\in \Net_{I,J}$ and a network $N(A)$ representing an element $A\in GL_n$ as explained above. We will concatenate $N$ and $N(A)$ according to the  following rule: in $N(A)$ reverse directions of all horizontal paths $i\to i $ for $i\in I$, changing every edge weight $w$ involved to $w^{-1}$, and then
 glue the left boundary of $N(A)$ to the boundary of $N$ in such a way that boundary vertices
 with the same label are glued to each other. Denote the resulting net by $N\circ N(A)$. Let $\wX_N$ and $\wX_{N\circ N(A)}$ be the signed boundary measurement matrices constructed according to the recipe outlined at the beginning of Section~\ref{indpoisson}.
 
\begin{lemma}
Matrices $\wX_{N\circ N(A)}$ and $\wX_NA$  are representatives of the same element in $G_k(n)$.
 \end{lemma}

\begin{proof}
To check that  $\wX_{N\circ N(A)}$ coincides with the result of the natural action of $GL_n$ on $G_k(n)$ induced by the right multiplication, it suffices to consider the case when $A$ is a diagonal or an elementary bidiagonal matrix, that is when $N(A)$ is one of the elementary networks in  Fig.~\ref{fig:factor}. If $N(A)$ is the first diagram in Fig.~\ref{fig:factor}, then the boundary measurements in $N'=N\circ N(A)$ are given by $M'(i, j) = d_i^{-1} M(i,j) d_j$,  $i\in I, j\in J $, where
 $M(i,j)$ are the boundary measurements in $N$. This is clearly consistent with the natural $GL_n$ action.  
 
 Now let $A=E_i^-(l)$ (the case $A=E_i^+(u)$ can be treated similarly). Then 
 \begin{equation}\label{M_A}
M'(i_p,j) = \begin{cases}  
M(i_p,j) + \delta_{i i_p} l M(i_{p-1},j) & \hskip -2.5cm\text{if $i,i-1 \in I$} \\
M(i_p,j) + \delta_{i-1, j} l M(i_p,j+1)  & \hskip -2.5cm\text{if $i,i-1 \in J$} \\
M(i_p,j) + \delta_{i i_p}\delta_{i-1,j} l& \hskip -2.5cm\text{if $i\in I$, $i-1 \in J$}\\
\dfrac{M(i-1,j)}{1 + l M(i-1,i)} 
(\delta_{i_p,i-1}\pm  (1-\delta_{i_p,i-1})l M(i_p,i))  &  \\ 
\hskip 1cm+(1-\delta_{i_p,i-1})M(i_p,j)                         
&\hskip -2.5cm\text{if $i\in J$, $i-1 \in I$}        
   \end{cases}
\end{equation}
 (in the last line above we used Lemma~\ref{recalc}). 
 
 Recall that for $j\in J$, an entry $\wX_{pj}$ of $\wX_N$ coincides with $M(i_p,j)$ 
 up to the  sign $(-1)^{s(p,j)}$. By the construction in Section~\ref{indpoisson}, if $ i,i-1 \in I $, then $(-1)^{s(p,j)}= - (-1)^{s(p-1,j)}$ for $p$ such that $i_p=i$ and for all $j\in J$. Then the first line of (\ref{M_A}) shows that 
 $\wX_{N'}= E_p^-(-l) \wX_N E_i^-(l)$, where the first elementary matrix is $k\times k$ and the second one is $n\times n$.  
 
If $ i,i-1 \in J $, then $(-1)^{s(p,i-1)}= (-1)^{s(p,i)}$ for all $p\in [1,k]$, and the second line of (\ref{M_A}) results in $\wX_{N'} = \wX_N E_i^-(l)$. The latter equality is also clearly valid for the case
 $ i\in I$, $i-1 \in J$. 
 
 Finally, if $ i\in J$, $i-1 \in I$, let $r\in [1,k]$ be the index such that $i-1=i_r$, and let
 $\wX_N^{(i)}$ denote the $i$th column of $\wX_N$. Then the $k\times k$ submatrix of $\wX_N E_i^-(l)$ formed by the columns indexed by $I$ is $C = \one_k + l \wX_N^{(i)} e_{r}^T$, where $e_r=(0,\dots,0,1,0,\dots,0)$ with the only nonzero element in the $r$th position. Therefore,
\[
 \tilde X = C^{-1} \wX_N E_i^-(l) = \left ( \one_k  - \frac{1}{ 1 + l \wX_{r i}} l \wX_N^{(i)} e_{r}^T \right ) \wX_N E_i^-(l) 
\]
is the representative of $[\wX_N E_i^-(l)]\in G_k(n)$ that has an identity matrix as its  $k\times k$ submatrix formed by columns indexed by $I$. We need to show that $\tilde X = \wX_{N'}$. First note that
 $\wX_{ri}=M(i_r, i)=M(i-1,i)$ and so, for all $j\in J$, 
\[
\tilde X_{rj} =   \frac{\wX_{rj}}{ 1 + l \wX_{r i}} =  \frac{(-1)^{s(r,j)}M(i-1,j)}{ 1 + l M(i-1,i)}= 
(-1)^{s(r,j)}M'(i-1,j)= (\wX_{N'})_{rj}.
\]
 If $p\ne r$, then
\begin{multline*}
 \tilde X_{pj} = \wX_{pj} -  \frac{l \wX_{rj} \wX_{pi} }{ 1 + l \wX_{r i}}\\
 = (-1)^{s(p,j)} \left( M(i_p, j) - (-1)^{s(p,i) + s(r,j)-s(p,j)} 
  \frac{l M(i_p, i) M(i-1,j) }{ 1 + l M(i-1,i)}\right ).
\end{multline*}

 Thus, to see that $\tilde X_{pj} = (\wX_{N'})_{rj}$, it is enough to check that the sign assignment that was used in Lemma~\ref{recalc} is consistent with the formula 
 \[
 - (-1)^{s(p,i) + s(r,j)-s(p,j)}.
 \]
  This can be done by direct inspection and is left to the reader as an exercise.
 \end{proof}
 
   Recall that if  $\H$ is a  Lie subgroup of a Poisson-Lie group $\G$, then a Poisson
structure on the homogeneous space $\H \backslash \G$ is called
{\it Poisson homogeneous\/} if the action map $\H \backslash \G
\times \G \to \H \backslash \G$ is Poisson.
 Now that we have established that the natural right action of $GL_n$ on $G_k(n)$ can be realized via the operation on networks described above,  
 Theorems~\ref{PSME},~\ref{PSGL} and~\ref{PSGr} immediately imply the following
 
 \begin{theorem}\label{poihomo}
 For any choice of parameters $\alpha, \beta$, the Grassmannian $G_k(n)$ equipped with the bracket $\P_{\alpha,\beta}$ is a Poisson homogeneous
 space for $GL_n$ equipped with the "matching" Poisson-Lie bracket~\eqref{Abrack12}.
\end{theorem} 

\begin{remark}{\rm
It is not difficult to deduce Theorem \ref{poihomo} from the general theory of Poisson homogeneous
spaces (see, e.g. \cite{r-sts}). The proof above illustrates capabilities of the network approach.}
\end{remark}

\vskip 1cm

\section{Compatibility with cluster algebra structure}
\label{compclust}

\subsection{Cluster algebras and compatible Poisson brackets}
\label{CA&PB}

First, we recall the basics of cluster algebras of geometric type. The definition that we present
below is not the most general one, see, e.g.,
\cite{FZ2, CAIII} for a detailed exposition.
 
The {\em coefficient group\/} $\PP$ is a free multiplicative abelian
group of a finite rank $m$ with generators $g_1,\dots, g_m$.
An {\em ambient field\/}  is
the field $\FFF$ of rational functions in $n$ independent variables with
coefficients in the field of fractions of the integer group ring
$\Z\PP=\Z[g_1^{\pm1},\dots,g_m^{\pm1}]$ (here we write
$x^{\pm1}$ instead of $x,x^{-1}$).

A {\em seed\/} (of {\em geometric type\/}) in $\FFF$ is a pair
$\Sigma=(\x,\widetilde{B})$,
where $\x=(x_1,\dots,x_n)$ is a transcendence basis of $\FFF$ over the field of
fractions of $\Z\PP$ and $\widetilde{B}$ is an $n\times(n+m)$ integer matrix
whose principal part $B$ (that is, the $n\times n$ submatrix formed by the
columns $1,\dots,n$) is skew-symmetric.

The $n$-tuple  $\x$ is called a {\em cluster\/}, and its elements
$x_1,\dots,x_n$ are called {\em cluster variables\/}. Denote
$x_{n+i}=g_i$ for $i\in [1,m]$. We say that
$\widetilde{\x}=(x_1,\dots,x_{n+m})$ is an {\em extended
cluster\/}, and $x_{n+1},\dots,x_{n+m}$ are {\em stable
variables\/}. It is convenient to think of $\FFF$ as
of the field of rational functions in $n+m$ independent variables
with rational coefficients. 

Given a seed as above, the {\em adjacent cluster\/} in direction $k\in [1,n]$
is defined by
$$
\x_k=(\x\setminus\{x_k\})\cup\{x'_k\},
$$
where the new cluster variable $x'_k$ is given by the {\em exchange relation}
\begin{equation}\label{exchange}
x_kx'_k=\prod_{\substack{1\le i\le n+m\\  b_{ki}>0}}x_i^{b_{ki}}+
       \prod_{\substack{1\le i\le n+m\\  b_{ki}<0}}x_i^{-b_{ki}};
\end{equation}
here, as usual, the product over the empty set is assumed to be
equal to~$1$.

We say that $\wB'$ is
obtained from $\wB$ by a {\em matrix mutation\/} in direction $k$
and
write $\wB'=\mu_k(\wB)$ 
 if
\[
b'_{ij}=\begin{cases}
         -b_{ij}, & \text{if $i=k$ or $j=k$;}\\
                 b_{ij}+\displaystyle\frac{|b_{ik}|b_{kj}+b_{ik}|b_{kj}|}2,
                                                  &\text{otherwise.}
        \end{cases}
\]

Given a seed $\Sigma=(\x,\widetilde{B})$, we say that a seed
$\Sigma'=(\x',\widetilde{B}')$ is {\em adjacent\/} to $\Sigma$ (in direction
$k$) if $\x'$ is adjacent to $\x$ in direction $k$ and $\widetilde{B}'=
\mu_k(\widetilde{B})$. Two seeds are {\em mutation equivalent\/} if they can
be connected by a sequence of pairwise adjacent seeds.

The {\em cluster
algebra\/} (of {\em geometric type\/}) $\A=\A(\widetilde{B})$
associated with $\Sigma$ is the $\Z\PP$-subalgebra of $\FFF$
generated by all cluster variables in all seeds mutation
equivalent to $\Sigma$.

Let $\Poi$ be a Poisson bracket on the ambient field $\FFF$
. We say that it is {\em compatible} with the cluster algebra $\A$ if, for any extended
cluster $\widetilde{\x}=(x_1,\dots,x_{n+m})$,  one has
$$\{x_i,x_j\}=\omega_{ij} x_ix_j\ ,$$
where $\omega_{ij}\in\Z$ are
constants for all $i,j\in[1,n+m]$. The matrix
$\Omega^{\widetilde \x}=(\omega_{ij})$ is called the {\it coefficient matrix\/}
of $\Poi$ (in the basis $\widetilde \x$); clearly, $\Omega^{\widetilde \x}$ is
skew-symmetric.

Consider, along with
cluster and stable variables $\wx$, another $(n+m)$-tuple of
rational functions denoted
$\tau=(\tau_1,\dots,\tau_{n+m})$ and defined by
\begin{equation}\label{eq:1.3}
\tau_j=x_j^{\varkappa_j}\prod_{k=1}^n x_k^{b_{jk}},
\end{equation}
where $\varkappa_j$ is an integer, $\varkappa_j=0$ for $1\ls j\ls n$.
We say that the entries $\tau_i,\, i\in[1,n+m]$ form
\emph{a $\tau$-cluster}. It is proved in \cite{GSV1}, Lemma 1.1, that 
$\varkappa_j, \ n+1 \leq j \leq n+m$ can be selected in such a way that
the transformation $\wx\mapsto\tau$ is  non-degenerate, provided 
$\rank \wB=n$.

Recall that a square matrix $A$ is {\it reducible\/} if there exists a
permutation matrix $P$ such that $PAP^T$ is a block-diagonal
matrix, and {\it irreducible\/} otherwise.
The following result is a particular case of Theorem 1.4 in \cite{GSV1}.

\begin{theorem}\label{thm:1.4}
Assume that $\rank \wB=n$ and the principal part of $\wB$ is irreducible.
Then a Poisson bracket is compatible with $\A(\wB)$ if and only if its coefficient matrix in the basis $\tau$ has the following property: its $n\times (n+m)$ submatrix formed by the first $n$ rows is proportional to
 $\wB$.
\end{theorem}

In \cite{GSV1} we constructed the cluster algebra $\A_{G^0_k(n)}$ on the open cell
$G^0_k(n)$ in the Grassmannian $G_k(n)$. This structure can be viewed as a restriction
of the cluster algebra of the coordinate ring of $G_k(n)$ described in \cite{Scott} using combinatorial
properties of Postnikov's construction.

\begin{figure}[ht]
\begin{center}
\includegraphics[height=4cm]{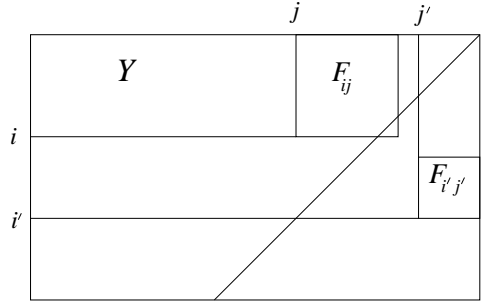}
\caption{To the definition of $F_{ij}$} \label{fviay}
\end{center}
\end{figure}

Let us briefly review the construction of \cite{GSV1}.
For every $(i,j)$-entry of the matrix $Y$  defined in (\ref{defY}), put
\begin{equation}
l(i,j)=\min (i-1, n-k-j)
\nonumber
\end{equation}
and
\begin{equation}\label{eq:3.16}
F_{ij}= Y_{[i-l(i,j),i]}^{[j, j+l(i,j)]}\ .
\end{equation}
Submatrices  of $Y$ whose determinants define functions $F_{ij}$ are depicted in Fig.~\ref{fviay}.

The initial extended cluster consists of
functions
\begin{multline}\label{eq:3.25}
  f_{ij} =(-1)^{(k-i)(l(i,j)- 1)} F_{ij}=
\frac{x_{([1,k]\setminus [i-l(i,j),\ i])
\cup [ j+k, j+l(i,j)+k]}}{x_{[1,k]}},\\
i\in [1,k],\quad j\in [1, m] \ ,
\end{multline}
where $x_I$ denote Pl\"ucker coordinates of the element of $G_k^0(n)$
represented by the matrix $\left [  \one_k \  Y  \right ]$.
Functions $f_{11},f_{21},\ldots, f_{k1},
f_{k2},\ldots, f_{km}$   serve as stable
coordinates.

\begin{figure}[ht]
\begin{center}
\includegraphics[height=3.5cm]{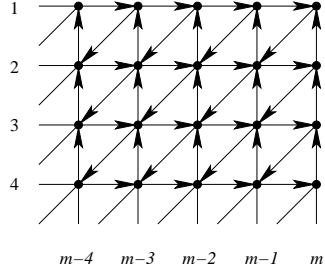}
\caption{Graph $\Gamma(\wB)$ corresponding to $G^0_k(n)$}\label{fig:graph}
\end{center}
\end{figure}

The entries of $\wB$ are all $0$ or $\pm 1$s. Thus it is convenient to describe
$\wB$ by a directed  graph $\Gamma(\wB)$. 
The
vertices of $\Gamma(\wB)$ correspond to all columns of $\wB$, and, since $\wB$ is rectangular, the corresponding edges are either  between the cluster variables or between a cluster variable and a stable variable. In our case, $\Gamma(\wB)$ is a directed graph
with vertices forming a rectangular $k\times m$ array and labeled
by pairs of integers $(i,j)$, and edges $(i,j) \to (i,j+1),\
(i+1,j) \to (i,j)$ and $(i,j) \to (i+1,j-1)$
(cf.~Fig.~\ref{fig:graph}).

The main goal of this section is to use networks in order to show that, in
fact, every Poisson structure in the 2-parameter family $\P_{\alpha,\beta}$ described in Theorem~\ref{PSGr} is compatible with $\A_{G^0_k(n)}$.

\subsection{Face weights and boundary measurements}
Let $N=(G,w)$ be a perfect planar network in a disk.
Graph $G$ divides the disk into a finite number of connected components called
\emph{faces}. The boundary of each face consists of edges of $G$ and, possibly, of several arcs bounding the disk. A face is called {\it bounded\/} if its boundary contains only edges of $G$ and {\it unbounded\/} otherwise. In this Section we additionally require that each edge of $G$ belongs to a path from a source to a sink. This is a technical condition that ensures that the two faces separated by an edge are distinct. Clearly, the edges that violate this condition do not influence the boundary measurement map and may be eliminated from the graph.  

Given a face $f$, we define its {\it face weight\/} $y_f$ as the Laurent monomial in edge weights $w_e$, $e\in E$, given by
\begin{equation}\label{faceweight}
y_f=\prod_{e\in\partial f}w_e^{\gamma_e},
\end{equation}
where $\gamma_e=1$ if the direction of $e$ is compatible with the counterclockwise orientation of the boundary $\partial f$ and $\gamma_e=-1$ otherwise. For example, the face weights for the network shown in Figure~\ref{fig:GrBMM} are $w_1w_2^{-1}w_3$, $w_3^{-1}w_5^{-1}w_8^{-1}$, $w_6w_7^{-1}w_8$, $w_1^{-1}w_4^{-1}w_6^{-1}$ for four unbounded faces and $w_2w_4w_5w_7$ for the only bounded face. 

Similarly to the space of edge weights $\CC_N$, we can define the {\it space of face weights\/} $\FF_N$; a point of this space is the graph $G$ as above with the faces weighted by real numbers obtained by specializing the variables $x_1,\dots,x_d$ in the expressions for $w_e$ to nonzero values and subsequent computation via~(\ref{faceweight}) (recall that by definition, $w_e$ do not vanish on $\EE_N$). By the Euler formula, the number of faces of $G$ equals $|E|-|V|+1$, and so $\FF_N$ is a semialgebraic subset in $\R^{|E|-|V|+1}$. In particular, the product of the face weights over all faces equals~$1$ identically, since each edge enters the boundaries of exactly two faces in opposite directions.

In general, the edge weights can not be restored from the face weights. However, the following proposition holds true.

\begin{lemma}\label{BMviay}
The weight $w_P$ of an arbitrary path $P$ between a source $b_i$ and a sink $b_j$ is a monomial in the face weights.
\end{lemma}

\begin{proof} Assume first that all the edges in $P$ are distinct. Extend $P$ to a cycle $C_P$ by adding the arc on the boundary of the disk between $b_j$ and $b_i$ in the counterclockwise direction. Then the product of the face weights over all faces lying inside $C_P$ equals the product of the conductivities over the edges of $P$, which is $w_P$. Similarly, for any simple cycle $C$ in $G$, its weight $w_C$ equals the product of the face weights over all faces lying inside $C$. It remains to use~(\ref{offcycle}) and to care that the cycles that are split off are simple. The latter can be guaranteed by the {\it loop erasure\/} procedure (see \cite{Loop}) that consists in traversing $P$ from $b_i$ to $b_j$ and splitting off the cycle that occurs when the first time the current edge of $P$ coincides with the edge traversed earlier.
\end{proof} 

It follows from Lemma~\ref{BMviay} that one can define boundary measurement maps $M_N^\FF\:\FF_N\to\Mat_{k,m}$ and $X_N^\FF\: \FF_N\to G_k(n)$ so that
$$
M_N^\FF\circ y=M_N,\qquad X_N^\FF\circ y=X_N,
$$
where $y\:\CC_N\to\FF_N$ is given by~(\ref{faceweight}).

\begin{figure}[ht]
\begin{center}
\includegraphics[height=4cm]{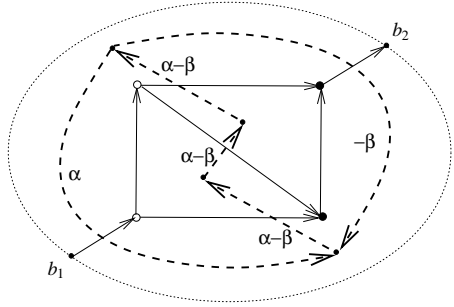}
\caption{Directed planar network and its directed dual} \label{fig:dualnet}
\end{center}
\end{figure}

Our next goal is to write down the 2-parametric family of Poisson structures induced on $\FF_N$ by the map $y$. Let $N=(G,w)$ be a perfect planar network. In this section it will be convenient to assume that boundary vertices are colored in {\it gray}.  Define the {\it directed dual network\/} $N^*=(G^*,w^*)$ as follows. Vertices of $G^*$ are the faces of $N$. Edges of $G^*$ correspond to the edges of $N$ with endpoints of different colors; note that there might be several edges between the same pair of vertices in $G^*$.
An edge $e^*$ of $G^*$ corresponding to $e$ is directed in such a way that the white endpoint of $e$ (if it exists) lies to the left of $e^*$ and the black endpoint of $e$ (if it exists) lies to the right of $e$. The weight $w^*(e^*)$ equals $\alpha-\beta$ if the endpoints of $e$ are white and black, $\alpha$ if the endpoints of $e$ are white and gray and $-\beta$ if the endpoints of $e$ are black and gray. An example of a directed planar network and its directed dual network is given in Fig.~\ref{fig:dualnet}.

\begin{lemma}\label{PSviay}
The $2$-parametric family $\{\cdot,\cdot\}_{\FF_N}$ is given by
\[
\{y_f,y_{f'}\}_{\FF_N}=\left(\sum_{e^*: f\to f'} w^*(e^*)-
\sum_{e^*: f'\to f} w^*(e^*)\right)y_fy_{f'}.
\]
\end{lemma}

\begin{proof}
Let $e=(u,v)$ be a directed edge. We say that the flag $(u,e)$ is {\it positive\/}, and the flag $(v,e)$ is {\it negative}. The color of a flag is defined as the color of the vertex participating in the flag.

Let $f$ and $f'$ be two faces of $N$. We say that a flag $(v,e)$ is {\it common\/} to $f$ and $f'$ if both $v$ and $e$ belong to $\partial f\cap\partial f'$. Clearly, the bracket  $\{y_f,y_{f'}\}_{\FF_N}$ can be calculated as the sum of the contributions of all flags common to $f$ and $f'$.

Assume that $(v,e)$ is a positive white flag common to $f$ and $f'$, see Fig.~\ref{fig:comflag}. Then $y_f=\dfrac{x_v^3}{x_v^2}\bar y_f$ and $y_{f'}=x_v^1x_v^2\bar y_{f'}$, where $x_v^i$ are the weights of flags involving $v$ and $\{x_v^i,\bar y_f\}_\RR=\{x_v^i,\bar y_{f'}\}_\RR=0$, see Section~\ref{PSSC}. Therefore, by~(\ref{6parw}), the contribution of $(v,e)$ equals $(\alpha_{12}-\alpha_{13}-\alpha_{23})y_fy_{f'}$, which by~(\ref{condw}) equals $-\alpha y_fy_{f'}$. 

Assume now that $(v,e)$ is a negative white flag common to $f$ and $f'$, see Fig.~\ref{fig:comflag}. In this case $y_f=\dfrac1{x_v^1x_v^3}\bar y_f$ and 
$y_{f'}=x_v^1x_v^2\bar y_{f'}$, so  the contribution of $(v,e)$ equals $(\alpha_{13}+\alpha_{23}-\alpha_{12})y_fy_{f'}=\alpha y_fy_{f'}$.

\begin{figure}[ht]
\begin{center}
\includegraphics[height=2.5cm]{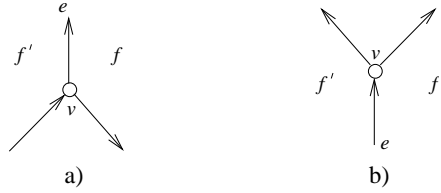}
\caption{Contribution of a white common flag: a) positive flag; b) negative flag} \label{fig:comflag}
\end{center}
\end{figure}

In a similar way one proves that the contribution of a positive black flag common to $f$ and $f'$ equals $-\beta y_fy_{f'}$, and the contribution of a negative black flag common to $f$ and $f'$ equals $\beta y_fy_{f'}$. Finally, the contributions of positive and negative gray flags are clearly equal to zero.

The statement of the lemma now follows from the definition of the directed dual network.
\end{proof}

\subsection{Compatibility theorem}
Now we can formulate the main  theorem of this section.

\begin{theorem}\label{thm:PoissComp}
Any Poisson structure in the two-parameter family $\P_{\alpha,\beta}$ is
compatible with the cluster algebra $\A_{G^0_k(n)}$.
\end{theorem}

\begin{proof} 
By Theorem~\ref{thm:1.4}, it is enough to choose an initial extended cluster and to compare the coefficient matrix of $\P_{\alpha,\beta}$ in the basis $\tau$ with the exchange matrix for this cluster. 

To find the coefficient matrix of $\P_{\alpha,\beta}$ we define a special network $N(k,m)$ with $k$ sources and $m$ sinks. The graph of $N(k,m)$ has $km+1$ faces. Each bounded face is a hexagon; all bounded faces together form a $(k-1)\times(m-1)$ parallelogram on the hexagonal lattice.
Edges of the hexagons are directed North, South-East and South-West.
 Each vertex of degree~2 on the left boundary of this parallelogram is connected to a source, and each vertex of degree~2 on the upper boundary of this parallelogram is connected to a sink. The remaining vertices of degree~2 on the lower and the right boundaries of the parallelogram are eliminated: two edges $u_1\to v\to u_2$ are replaced by one edge $u_1\to u_2$ and the intermediate vertex $v$ is deleted. The sources are labelled counterclockwise from~1 to $k$, sinks from $k+1$ to $n$. The faces are labelled by pairs $(ij)$ such that $i\in [1,k]$, $j\in [k+1, n]$. The unbounded faces are labelled $(1j)$ and $(in)$ except for one face, which is not labelled at all. The network $N(k,m)$ is defined by assigning a face weight $y_{ij}$ to each labelled face $(ij)$. Consequently, the directed dual network $N^*(k,m)$ forms the dual triangular lattice. All edges of $N^*(k,m)$ incident to the vertices $(ij)$, $i\in [2,k]$, $j\in [k+1, n-1]$ are of weight $\alpha-\beta$.
 
 For an example of the construction for $k=3$, $m=4$ see Fig.~\ref{fig:postnet}. The edge weights of the dual network that are not shown explicitly are equal to $\alpha-\beta$. 

\begin{figure}[ht]
\begin{center}
\includegraphics[height=5cm]{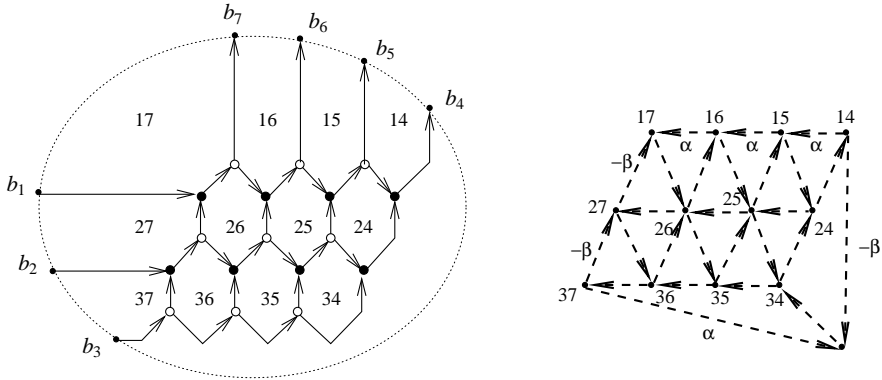}
\caption{The graph of $N(3,4)$ and its directed dual $N^*(3,4)$} \label{fig:postnet}
\end{center}
\end{figure}

As the first step of the proof, we express the cluster variables $f_{ij}$ via the face weights of $N(k,m)$.

\begin{lemma} \label{pmfviay}
For any $i\in [1,k]$, $j\in [1,m]$, one has
\begin{equation}\label{eq:pmfviay}
f_{ij}=\pm\prod_{p=1}^i\prod_{q=j+k}^n y_{pq}^{1+\min\{i-p,q-j-k\}}.
\end{equation}
\end{lemma}

\begin{proof}
Since the graph of $N(k,m)$ is acyclic, one can use the Lindstr\"om lemma \cite{Lindstrom} to calculate the minors of the boundary measurement matrix. Assume first that $i+j\le n+1-k$. By~(\ref{eq:3.16}) and the Lindstr\"om lemma, $f_{ij}$ equals to the sum of the products of the path weights for all $i$-tuples of nonintersecting paths between the sources $b_1,b_2,\dots,b_i$ and the sinks $b_{j+k}, b_{j+k+1}, \dots, b_{j+k+i-1}$, each product being taken with a certain sign. Note that there exists a unique path from $b_1$ to $b_{i+j+k-1}$ in $N(k,m)$. After this path is chosen, there remains a unique path from $b_2$ to $b_{i+j+k-2}$, and so on. Moreover, a path from $b_1$ to any sink $b_p$ with $p<i+j+k-1$ cuts off the sink $b_{i+j+k-1}$. Therefore, there exists exactly one $i$-tuple of paths as required. Relation~(\ref{eq:pmfviay}) 
for $i+j\le n+1-k$ now follows from the proof of Lemma~\ref{BMviay}.

The case $i+j\ge n+1-k$ is similar to the above one and relays on the uniqueness of an $(n+1-k-j)$-tuple of nonintersecting paths between the sources $b_{i+j+k-n},
 b_{i+j+k-n+1},\dots, b_i$ and the sinks $b_{j+k}, b_{j+k+1}, \dots, b_n$. Here we start with the unique path from $b_i$ to $b_n$, then choose the unique remaining path between $b_{i-1}$ and $b_{n-1}$, and so on. 
\end{proof}

The next step is the calculation of $\tau$-coordinates for the initial extended cluster~(\ref{eq:3.25}) via face weights.

\begin{lemma}\label{tauviay1}
The $\tau$-coordinates corresponding to the cluster variables of the initial extended cluster~\eqref{eq:3.25} are given by
\begin{equation}\label{eq:tauviay1}
\tau_{ij}=\pm y_{i+1,j+k-1}, \qquad i\in [1,k-1],\; j\in [2,m].
\end{equation}
\end{lemma}

\begin{proof}
Combining the definition of $\tau$-coordinates via~(\ref{eq:1.3}) and the description of the exchange matrix in the basis $\{f_{ij}\}$ provided by Fig.~\ref{fig:graph}, we conclude that
\[
\tau_{ij}=\frac{f_{i+1,j-1}f_{i,j+1}f_{i-1,j}}{f_{i,j-1}f_{i+1,j}f_{i-1,j+1}}=
\frac{\dfrac{f_{i+1,j-1}}{f_{ij}}\cdot\dfrac{f_{ij}}{f_{i-1,j+1}}}{\dfrac{f_{i+1,j}}{f_{i,j+1}}\cdot\dfrac{f_{i,j-1}}{f_{i-1,j}}}, \qquad
i\in [1,k-1],\; j\in [2,m];
\]
here we assume that $f_{ij}=1$ if $i=0$ or $j=m+1$. Next, by~(\ref{eq:pmfviay}),
\[
\frac{f_{ij}}{f_{i-1,j+1}}=\pm\prod_{p=1}^i\prod_{q=j+k}^n y_{pq}
\]
for $i\in [1,k-1]$, $j\in [2,m]$, and the result follows.
\end{proof}

\begin{remark} The operation that expresses $\tau$-coordinates in terms of face weights is an analog of the twist studied in \cite{BFZ,FZBruhat}.
\end{remark}

The expressions for the $\tau$-coordinates that correspond to the stable variables of the 
initial extended cluster~(\ref{eq:3.25}) are somewhat cumbersome. Recall that by~(\ref{eq:1.3}), each stable variable enters the expression for the corresponding $\tau$-coordinate with some integer exponent $\varkappa_{ij}$. Let us denote
\[
\tau^*_{ij}=\tau_{ij}f_{ij}^{-\varkappa_{ij}}
\]
for $i=k$, $j\in [1,m]$ and $i\in [1,k]$, $j=1$.

\begin{lemma}\label{tauviay2}
The $\tau$-coordinates corresponding to the stable variables of the initial extended cluster~\eqref{eq:3.25} are given by
\begin{equation}\label{eq:tauviay2}
\tau^*_{ij}=\begin{cases}
\pm\displaystyle\prod_{p=1}^{k-1}\;\prod_{q=j+k}^{\min\{n,j+k-p-1\}}y_{pq} \qquad &\text{for $i=k$, $j\in [2,m-1]$},\\
\pm\displaystyle\prod_{q=k+2}^n\;\prod_{p=i}^{\min\{n, i-q+k+2\}}y_{pq} \qquad &\text{for $i\in [2,k-1]$, $j=1$},\\
\pm\displaystyle\prod_{p=1}^{k-1}y_{pn} \qquad &\text{for $i=k$, $j=m$},\\
\pm\displaystyle\prod_{q=k+2}^{n}y_{1q} \qquad &\text{for $i=1$, $j=1$},\\
\pm\displaystyle\prod_{p=1}^k\prod_{q=k+1}^n y_{pq}^{-\min\{k-p,q-k-1\}} \qquad &\text{for $i=k$, $j=1$}.
           \end{cases}
\end{equation}
\end{lemma}

\begin{proof}
It suffices to note that
\[
\tau^*_{ij}=\begin{cases}
          f_{k-1,j}/f_{k-1,j+1} \qquad &\text{for $i=k$, $j\in [2,m-1]$},\\
          f_{i,2}/f_{i-1,2} \qquad &\text{for $i\in [2,k-1]$, $j=1$},\\
          f_{k-1,m} \qquad &\text{for $i=k$, $j=m$},\\
          f_{12} \qquad &\text{for $i=1$, $j=1$},\\
          1/f_{k-1,2} \qquad &\text{for $i=k$, $j=k+1$},
           \end{cases}
\]
and to apply Lemma~\ref{pmfviay}.
\end{proof}

To conclude the proof of the theorem we build the network $\Gamma_{\alpha,\beta}$ representing the family of Poisson brackets $\P_{\alpha,\beta}$ in the basis $\tau$. The vertices of $\Gamma_{\alpha,\beta}$ are $\tau$-coordinates $\tau_{ij}$, $i\in [1,n]$, $j\in [1,m]$. An edge from $\tau_{ij}$ to $\tau_{pq}$ with weight $c$ means that $\{\tau_{ij},\tau_{pq}\}_{\FF_{N(k,m)}}=c\tau_{ij}\tau_{pq}$. 

It follows immediately from Lemmas~\ref{PSviay} and~\ref{tauviay1} that the induced subnetwork $\Gamma^0_{\alpha,\beta}$ of $\Gamma_{\alpha,\beta}$ spanned by the vertices $\tau_{ij}$, $i\in [1,k-1]$, $j\in [2,m]$ is isomorphic to the subnetwork $N_0^*(k,m)$ of $N^*(k,m)$ spanned by the vertices $(pq)$, $p\in [2,k]$, $q\in [k+1,n-1]$. We thus get an isomorphism between $\Gamma^0_{\alpha,\beta}$ and the corresponding induced subgraph of  $\Gamma(\wB)$, see Fig.~\ref{fig:iso}. Under this isomorphism each edge of weight $\alpha-\beta$ is mapped to an edge of weight~1.

\begin{figure}[ht]
\begin{center}
\includegraphics[height=4cm]{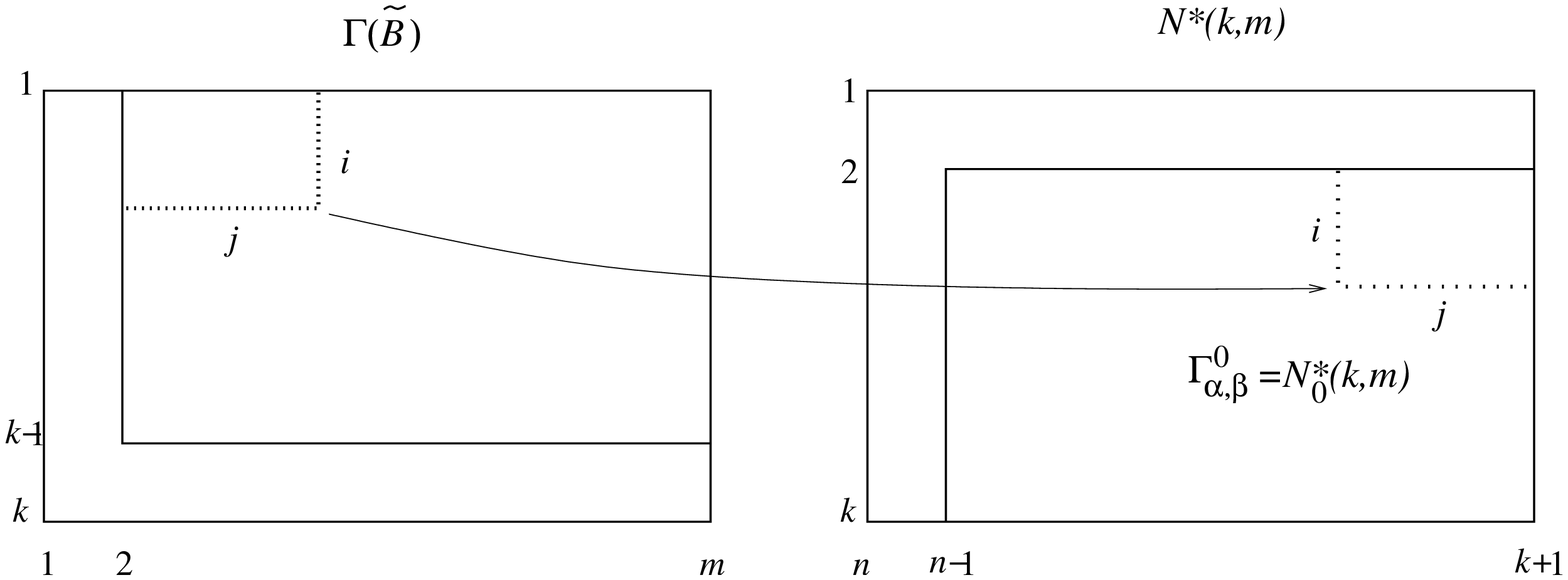}
\caption{The isomorphisms between induced subgraphs of  $\Gamma_{\alpha,\beta}$ and $\Gamma(\wB)$} \label{fig:iso}
\end{center}
\end{figure}

For the remaining vertices of $\Gamma_{\alpha,\beta}$, that is, $\tau_{ij}$, $i\in [1,k]$, $j=1$ or $i=k$, $j\in [1,m]$, we have to check the edges connecting them to the vertices of $\Gamma^0_{\alpha,\beta}$. 
Consider for example the case $i=k$, $j\in [2,m-1]$. Then, by Lemma~\ref{tauviay2}, the bracket $\{\tau^*_{kj},\tau_{pq}\}_{\FF_{N(k,m)}}$ 
for $\tau_{pq}\in \Gamma^0_{\alpha,\beta}$ 
is defined by the edges incident to the vertex subset of $N^*(k,m)$  corresponding to the factors in the right hand side of the first formula in~(\ref{eq:tauviay2}), see Fig.~\ref{fig:contract}. Clearly, any vertex 
$(p,q)\in N_0^*(k,m)$ other than $(k,j+k)$ and $(k,j+k-1)$ is connected to this subset by an even number of edges (more exactly, 0, 2, 4 or 6), all of them of weight $\alpha-\beta$. Since exactly half of the edges are directed to $(pq)$, the bracket between $y_{pq}$ and  $\tau^*_{kj}$ vanishes by Lemma~\ref{PSviay}. The remaining two edges connecting the contracted set to the vertices $(k,j+k)$ and $(k,j+k-1)$ correspond to the two edges connecting the vertex $(k,j)$ to $(k-1,j)$ and $(k-1,j+1)$ in $\Gamma(\wB)$, see  Fig.~\ref{fig:contract}. 

\begin{figure}[ht]
\begin{center}
\includegraphics[height=4cm]{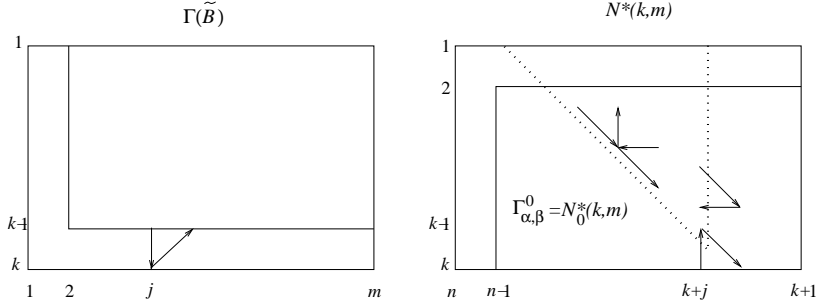}
\caption{The isomorphisms in the neighborhood of $\tau_{kj}$} \label{fig:contract}
\end{center}
\end{figure}

Finally, factor $f_{kj}$ commutes with all $\tau_{pq}\in\Gamma^0_{\alpha,\beta}$. 
Indeed, given an edge between $\tau_{pq}\in \Gamma^0_{\alpha,\beta}$ and $y_{rs}$, define its degree as the exponent of $y_{rs}$ in expression~(\ref{eq:pmfviay}) for $f_{kj}$. We have to prove that the sum of the degrees of all edges entering $\tau_{pq}$ equals the sum of the degrees of all edges leaving $\tau_{pq}$. This fact can be proved by analyzing all possible configurations of edges, see Fig.~\ref{fig:baldegree} presenting these configurations up to reflection.

\begin{figure}[ht]
\begin{center}
\includegraphics[height=3cm]{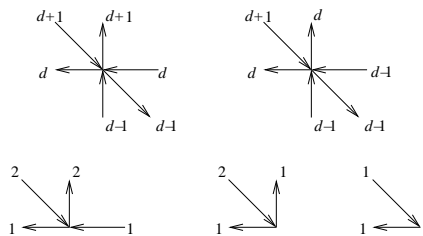}
\caption{Edges incident to $\tau_{pq}$ and their degrees} \label{fig:baldegree}
\end{center}
\end{figure}

Other cases listed in Lemma~\ref{tauviay2} are treated in the same way.
\end{proof}

\section{Acknowledgments}

We wish to express gratitude to A. Postnikov who explained to us the details
of his construction and to S. Fomin for stimulating discussions. M. G. was supported in part by NSF Grant DMS \#0400484. M. S.  was supported in part by NSF Grants DMS \#0401178 and PHY\#0555346.  Authors were
also supported by the BSF Grant \# 2002375. A. V. is grateful to Centre Interfacultaire Bernoulli at \'Ecole Polytechnique F\'ed\'erale de Lausanne for hospitality during his Spring 2008 visit.    We are also
grateful to the Institute of Advanced Studies of the Hebrew University of Jerusalem for the invitation to the Midrasha
Mathematicae  in May 2008, during which this paper was completed.


\begin{thebibliography}{O}

\bibitem [B] {Brenti} F. Brenti, \textit{Combinatorics and total positivity}, J. Combin. Theory Ser. A
{\bf 71} (1995), 175--218.

\bibitem [BFZ1] {BFZ}  A.~Berenstein, S.~Fomin, and A.~Zelevinsky,
\textit{Parametrizations of canonical bases and totally positive
matrices}. 
Adv. Math. \textbf{122} (1996), 49--149.

\bibitem [BFZ2] {CAIII}  A.~Berenstein, S.~Fomin, and A.~Zelevinsky,
\textit{Cluster algebras. III. Upper bounds and double Bruhat cells}. 
Duke Math. J. \textbf{126} (2005), 1--52.

\bibitem[Fa]{Fallat} S.~Fallat, \textit{ Bidiagonal factorizations of totally nonnegative matrices},
Amer. Math. Monthly, {\bf 108} (2001), 697--712.

\bibitem[Fo]{Loop} S.~Fomin, \textit{ Loop-erased walks and total positivity}, Trans. Amer. Math. Soc.
{\bf 353} (2001), 3563--3583.

\bibitem[FG1] {FayGekh1} L. Faybusovich,  M. I. Gekhtman,
\textit{Elementary Toda orbits and integrable lattices}, J. Math. Phys. {\bf 41} (2000),
2905--2921.

\bibitem[FG2] {FayGekh2} L. Faybusovich,  M. I. Gekhtman,\textit{
Poisson brackets on rational functions and multi-Hamiltonian
structure for integrable lattices}, Phys. Lett. A {\bf 272} (2000),
236--244.


\bibitem [FZ1] {FZBruhat}  S.~Fomin and A.~Zelevinsky, 
\textit{Double Bruhat cells and total positivity}. 
J. Amer. Math. Soc. \textbf{12}
(1999), 335--380.

\bibitem [FZ2] {FZ_Intel}  S.~Fomin and A.~Zelevinsky, 
\textit{Total Positivity: tests and parametrizations}.,
Math. Inteligencer. \textbf{22} 
(2000), 23--33.

\bibitem [FZ3] {FZ2}  S.~Fomin and A.~Zelevinsky, 
\textit{Cluster algebras.I. Foundations}. 
J. Amer. Math. Soc. \textbf{15} (2002), 497--529.



\bibitem [GSV1] {GSV1}  M.~Gekhtman, M.~Shapiro, and A.~Vainshtein,
\textit{Cluster algebras and Poisson geometry}.  
Mosc. Math. J. \textbf{3} (2003), 899--934.

\bibitem [GSV2] {GSV2} M.~Gekhtman, M.~Shapiro, and A.~Vainshtein,
\textit{Cluster algebras and Weil-Petersson forms}.
 Duke Math. J. \textbf{127} (2005), 291--311.


\bibitem[GSV3]{GSV3}  M.~Gekhtman, M.~Shapiro, and A.~Vainshtein, \textit{Poisson geometry of directed networks in an annulus}, arXiv: 0901.0020.

\bibitem[GSV4]{GSV4}   M.~Gekhtman, M.~Shapiro, and A.~Vainshtein,  \textit{  B\"acklund-Darboux transformations for Toda flows from cluster algebra perspective}, in preparation.

\bibitem[GSV5]{GSV5}  M.~Gekhtman, M.~Shapiro, and A.~Vainshtein, \textit{  Inverse problem for 1-1 networks in an annulus}, in preparation.

\bibitem [GrSh] {GrSh} B.~Gr\" unbaum and G.~Shephard,
\textit{Rotation and winding numbers for planar polygons and curves}.  
Trans. Amer. Math. Soc.  \textbf{322}  (1990), 169--187.

\bibitem[HKKR]{Reshetikhin&Co} T.~ Hoffmann, J. ~Kellendonk, N.~Kutz and N.~Reshetikhin \textit{Factorization dynamics and   Coxeter-Toda lattices}, Comm. Mat. Phys. {\bf 212} (2000), 297--321.


\bibitem[KM]{KarlinMacGregor} S. Karlin, J. McGregor,\textit{Coincidence Probabilities}, Pacific J. Math.
{\bf 9} (1959), 1141--1164.

\bibitem[LP]{Pavlo} T. Lam, P. Pylyavskyy, \textit{Total positivity in loop groups I: whirls and curls}, arXiv:0812.0840.

\bibitem[Li]{Lindstrom} B. Lindstr\"om, \textit{ On the vector representations of induced matroids}, Bull. London Math. Soc. {\bf 5} (1973), 85--90.

\bibitem[Lu]{Lusztig} G. Lusztig, \textit{ Total positivity in reductive groups}. Lie theory and geometry, 531--568, Progr. Math. {\bf 123}, Birkh\"auser, Boston, 1994.

\bibitem[M]{moser} J. Moser, \textit{Finitely many mass points on the line under the influence of the exponential potential - an integrable system}. Dynamical systems, theory and applications, 467--497,
Lecture Notes in Physics {\bf 38}, Springer, Berlin, 1975.

\bibitem[P]{Postnikov} A. Postnikov, \textit{ Total positivity, Grassmannians and networks}, arXiv: math/0609764.



\bibitem [R] {R2}  N.~Reshetikhin \textit{Integrability of
characteristic Hamiltonian systems on simple Lie groups with
standard Poisson Lie structure}, Comm. Mat. Phys. {\bf 242} (2003), 1--29.

\bibitem [ReST] {r-sts}  A.~Reyman and M.~Semenov-Tian-Shansky
\textit{Group-theoretical methods in the theory of
finite-dimensional integrable systems} Encyclopaedia of
Mathematical Sciences, vol.16, Springer--Verlag, Berlin, 1994 pp.
116--225.

\bibitem[S]{Scott} J. Scott, \textit{Grassmannians and cluster algebras}, Proc. London Math. Soc.
{\bf 92} (2006), 345--380.


\bibitem [T]{Talaska} K.~Talaska, \textit{A formula for Pl\"ucker coordinates associated with a planar network}, Int. Math. Res. Not. {2008}, 19pp.

\end{thebibliography}
 \end{document}